\documentclass[11pt]{article}
\usepackage{amsmath,amsfonts,amssymb,amsthm,mathrsfs}
\usepackage[a4paper,vmargin={3.5cm,3.5cm},hmargin={2.5cm,2.5cm}]{geometry}
\usepackage[font=sf, labelfont={sf,bf}, margin=1cm]{caption}
\usepackage{graphicx,graphics}
\usepackage{epsfig}
\usepackage{latexsym}
\usepackage[applemac]{inputenc}
\usepackage{ae,aecompl}
\usepackage[english]{babel}
\usepackage[colorlinks=true]{hyperref}
\usepackage{pstricks}
\usepackage{enumerate}
\usepackage{bbm}

\date{}

\newtheorem{theorem}{Theorem}[]
\newtheorem{remark}{Remark}[]

\newtheorem{proposition}[theorem]{Proposition}
\newtheorem{lemma}[theorem]{Lemma}
\newtheorem{corollary}[theorem]{Corollary}
\newtheorem{open}[]{Open question}

\def\llbracket{[\hspace{-.10em} [ }
\def\rrbracket{ ] \hspace{-.10em}]}

\newcommand{\eps}{\varepsilon}

\title{\bf \textsc{Random trees constructed by aggregation} \\
(Arbres aléatoires construits par agrégation)}

\author{\text{Nicolas Curien}\thanks{Universit\'e Paris--Sud, E-mail: nicolas.curien@gmail.com}  \ \ \& \text{B\'en\'edicte  Haas}\thanks{ Universit\'e Paris--Dauphine, E-mail: haas@ceremade.dauphine.fr} }

\begin{document}
\maketitle

We study a general procedure that builds random $\mathbb R$-trees by gluing recursively a new branch on a uniform point of the pre-existing tree. The aim of this paper is to see how the asymptotic behavior of the sequence of lengths of branches influences some geometric properties of the limiting tree, such as compactness and Hausdorff dimension. In particular, when the sequence of lengths of branches behaves roughly like $n^{-\alpha}$ for some $\alpha \in (0,1]$, we show that the limiting tree is a compact random tree of Hausdorff dimension $\alpha^{-1}$. This encompasses the famous construction of the Brownian tree of Aldous. When $\alpha >1$, the limiting tree is thinner and its Hausdorff dimension is always 1. In that case, we show that $ \alpha^{-1}$ corresponds to the dimension of the set of leaves of the tree.

\medskip

{\center{\textbf{Résumé}}

\medskip

Nous nous intéressons à une procédure générale de construction d'arbres réels aléatoires par collages successifs de nouvelles branches. A chaque étape, la nouvelle branche est collée en un point choisi uniformément sur l'arbre pré-existant. Notre objectif principal est de comprendre comment le comportement asymptotique de la suite des longueurs de branches influence certaines propriétés géométriques de l'arbre, telles que la compacité ou la dimension de Hausdorff. Nous montrons en particulier que lorsque la suite de longueurs de branches se comporte en $n^{-\alpha}$, avec $\alpha \in (0,1]$ fixé, l'arbre limite est compact, de dimension de Hausdorff $\alpha^{-1}$.  A titre d'exemple, ceci englobe une construction bien connue de l'arbre brownien d'Aldous. Lorsque $\alpha>1$, l'arbre limite est plus fin et de dimension de Hausdorff 1. Dans ce cas, nous montrons que $\alpha^{-1}$ correspond à la dimension de l'ensemble des feuilles de l'arbre.}

\section*{Introduction}
Consider a sequence of closed segments or ``branches'' of lengths  $ \ a_{1}, a_{2}, a_{3},... > 0$ and let
$$A_i=a_1+\ldots +a_i, \quad i \geq 1$$
denote the partial sums of their lengths.
We construct a sequence of random trees $ (\mathcal{T}_{n})_{n \geq 1}$ by starting with the tree $ \mathcal{T}_{1}$ made of the single branch of length $a_{1}$ and then recursively gluing the branch of length $a_{i}$ on a point uniformly distributed (for the length measure) on $ \mathcal{T}_{i-1}$. Let $ \mathcal{T}$ be the completion of the increasing union of the $ \mathcal{T}_{n}$ which is thus a random complete  continuous tree.  The aim of this paper is to discuss some geometric properties of this tree. Our first result shows that even if the series $\sum a_{i}$ is divergent, provided that  the sequence  $\mathbf a=(a_i)_{i\geq 1}$ is sufficiently well-behaved, the tree $ \mathcal{T}$ is a compact random tree with a fractal behavior.
\begin{theorem}[Case $\alpha \leq 1$] \label{thm:main} Suppose that there exists $\alpha \in (0,1]$ such that 
 \begin{eqnarray*} a_{i} \leq i^{-\alpha +\circ(1)} \qquad \mbox{and} \qquad   A_i= i^{1-\alpha + \circ(1)} \quad \mbox{ as }i \to \infty.  \end{eqnarray*} Then $ \mathcal{T}$ is almost surely a \emph{compact} real tree of Hausdorff dimension $ \alpha^{-1}$.
\end{theorem}

We actually get  more complete results. On the one hand, the tree $\mathcal T$ is  compact and has a Hausdorff dimension at most $\alpha^{-1}$ as soon as  $a_i \leq i^{-\alpha+\circ(1)}$ for some $\alpha \in (0,1]$ (Proposition \ref{prop:compact}). On the other hand, its Hausdorff dimension  is at least $\alpha^{-1}$ as soon as $A_i \geq i^{1-\alpha +\circ(1)}$  for some $\alpha \in (0,1]$ (Proposition \ref{prop:lowerbound} -- this result actually holds under a mild additional assumption that will be discussed in the core of the paper). Let us also mentioned that in a recent paper  \cite{ADGO14}, Amini et al.~considered the same aggregation model and obtained a necessary and sufficient condition for $ \mathcal{T}$ to be bounded in the particular case when $ \mathbf{a}$ is decreasing, see the discussion in Section \ref{sec:bounded}.

Theorem \ref{thm:main} encompasses the famous line-breaking construction of the Brownian continuum random tree (CRT) of Aldous. Specifically, if the sequence $\mathbf a$ is the random sequence of lengths given by the intervals in a Poisson process on $ \mathbb{R}_{+}$ with intensity $t \,\mathrm{d}t$, then Aldous proved \cite{Ald91a} that $ \mathcal{T}$ is compact and of Hausdorff dimension $2$ (this was the initial definition of the Brownian CRT). Yet, it is a simple exercise to see that such sequences almost surely satisfy the assumptions of our theorem for $\alpha = 1/2$. More generally, random trees built from a sequence of branches given by the intervals of a Poisson process of intensity $t ^\beta \mathrm{d}t$ on $ \mathbb{R}_{+}$ with $ \beta >0$ satisfy our assumptions with $\alpha = \beta/(\beta+1)$. Typically, in these examples, the sequence $\mathbf a$ is not monotonic.

\medskip

When the series $\sum a_{i}$ is convergent the situation may seem easier. In such cases, it should be intuitive that the limiting tree is compact  and of Hausdorff dimension $1$. We will see that this is true regardless of the mechanism used to glue the branches together (Proposition \ref{prop:finitelength}). But we can go further: when the asymptotic behavior of the sequence $\mathbf a$ is sufficiently regular, the set of leaves of $ \mathcal{T}$ exhibits an interesting fractal behavior similar to Theorem \ref{thm:main}. We recall that the leaves of a continuous tree $ \mathcal{T}$ are the points $ x$ such that $ \mathcal{T} \backslash \{x\}$ stays connected. 
\begin{theorem}[Case $\alpha > 1$]  \label{thm:maina>1} Suppose that there exists $\alpha >1 $ such that 
 \begin{eqnarray*} a_{i} \leq i^{-\alpha +\circ(1)} \qquad \mbox{and} \qquad  a_{i} + a_{i+1} + ... +a_{2i}= i^{1-\alpha + \circ(1)}  \quad \mbox{ as }i \to \infty.  \end{eqnarray*} Then the set of leaves of $ \mathcal{T}$ is almost  surely of Hausdorff dimension $ \alpha^{-1}$.
 \end{theorem}

We can decompose the tree $\mathcal T$ into its set of leaves $\mathsf{Leaves}(\mathcal T)$ and its skeleton $\mathcal T \backslash \mathsf{Leaves}(\mathcal T)$. Since the skeleton is a countable union of segments, its Hausdorff dimension is $1$ and so $\dim_{\mathrm H}(\mathcal T) = \linebreak  1 \vee \dim_{\mathrm H}\left(\mathsf{Leaves}(\mathcal T)\right)$.  Theorem \ref{thm:main} and  Theorem \ref{thm:maina>1} thus imply that when $a_i=i^{-\alpha}$ for some $\alpha \in (0,\infty)$, the tree $\mathcal T$ is compact and
$$\dim_{\mathrm H}\left(\mathsf{Leaves}(\mathcal T) \right)=\alpha^{-1}$$ 
almost surely.  When $\alpha=1$, the Hausdorff dimension of the leaves of $\mathcal T$ is not explicitly given in these theorems, but will be calculated further in the text. 

\medskip

\noindent \textsc{A toy-model for DLA.} Apart from the abundant random tree literature and the initial definition of the Brownian CRT by Aldous, a motivation for considering the above line-breaking construction is that it can be seen as a toy model of external diffusion limited aggregation (DLA). Recall that in the standard DLA model, say on $ \mathbb{Z}^2$, a subset $ \mathcal{A}_{n}$ is grown by recursively adding at each time a site on the boundary of $ \mathcal{A}_{n}$ according to the harmonic measure from infinity. It still remains a challenging open problem to understand the growth of $ \mathcal{A}_{n}$, see \cite{BPP97,Sch06}. In our model the particles are now branches of varying size (we do not rescale the aggregate) and harmonic measure seen from infinity is replaced by uniform measure on the structure at time $n$. Our Theorem \ref{thm:main} can thus be interpreted as the fact that in this case the DLA aggregate does not grow arms towards infinity, and identifies its fractal dimension.  

\medskip

The article \cite{HExplosion16} completes the previous results by studying cases where the tree $\mathcal T$ is obviously unbounded. Assuming that $(a_i)$ is regularly varying with a positive index, it describes the asymptotic behavior of the height of $\mathcal T_n$ and of the subtrees of $\mathcal T_n$ spanned by $\ell$ points picked uniformly and independently in $\mathcal T_n$, for all $\ell \in \mathbb N$. In another direction, Sénizergues  \cite{Delphin16+} extends our results to random metric spaces constructed by aggregation of $d$-dimensional spheres or more general independent random measured metric spaces, with gluing rules that depend both on the diameters and the measures of the metric spaces.  He shows an unexpected and intriguing Hausdorff dimension. Last we mention \cite{GH14} for a recent construction of the so-called stable trees \emph{via} an aggregation procedure that generalizes the line-breaking construction of the Brownian CRT, but that does not exactly fall in our setup.

\medskip

We finish this introduction by giving some elements of the proofs of our main results. In that aim, introduce the  quantity 
$$  \mathsf{H}( \mathbf{a}) \quad := \quad \sum_{i=1}^\infty \frac{a_{i}^2}{A_{i}}.$$ When the sequence $\mathbf a$ is bounded, we will see (Theorem \ref{thm:measures}) that condition $ \mathsf{H}( \mathbf{a}) < \infty$ is equivalent to the convergence of the normalized length measure $\mu_{n}$ on $ \mathcal{T}_{n}$ towards a limiting random probability $\mu$ on $ \mathcal{T}$. For connoisseurs, the latter is equivalent to the convergence of  $ (\mathcal{T}_{n},\mu_{n})$ to $( \mathcal{T}, \mu)$ in the Gromov--Prokhorov sense.  In particular, condition $ \mathsf{H}( \mathbf{a})< \infty$ ensures that the height of a  ``typical'' point of $ \mathcal{T}$ (i.e.~sampled according to $\mu$) is bounded. However it does not prevent $ \mathcal{T}$ from having very thin tentacles making it unbounded. 

Under the hypotheses of Theorem \ref{thm:main}, this phenomenon cannot happen thanks to an approximate scale invariance of the process. Roughly speaking, we prove that when $a_i \leq i^{-\alpha+\circ(1)}$, the subtree descending from the $i$th branch is a random tree built by an aggregation process which is similar to the construction of the original tree except that it is scaled by a factor at most $i^{-\alpha+\circ(1)}$. This gives the first hint that the fractal dimension of $ \mathcal{T}$ is at most $\alpha^{-1}$. On the other hand, when $A_i \geq i^{1-\alpha+\circ(1)}$ and $ \mathsf{H}( \mathbf{a})<\infty$, the lower bound on the dimension is obtained using Frostman's theory by constructing a (random) measure nicely spread  on $ \mathcal{T}$. This role will be played by the limiting measure $\mu$. To estimate the $\mu$-measure of typical balls of radius $r>0$ in $\mathcal{T}$ (Lemma  \ref{lem:majoDn}) we will compute the distribution of the distance of two typical points picked independently at random according to $\mu$ in $\mathcal T$, a.k.a.~the two-point function (Lemma \ref{lem:Dinfini}).

Under the hypotheses of Theorem \ref{thm:maina>1}, the upper bound of the dimension of the set of leaves is  even true in a deterministic setting (Proposition \ref{prop:finitelength}), as well as the compactness, and is obtained by exhibiting appropriate coverings. The lower bound of the dimension is again obtained via Frostman's theory. A difficulty in this case is that the random measure $\mu$ is equal to the normalized length measure on $ \mathcal{T}$ (recall that the total length of $ \mathcal{T}$ is finite in this case). Hence, $\mu$ is supported by the skeleton of the tree, and not by the leaves. This forces us to introduce another random measure supported by the leaves of $ \mathcal{T}$ which captures its fractal behavior.  This is done in the last section which is maybe the most technical part of this work. \medskip 

\noindent \textbf{Acknowledgments:} We thank the organizers and the participants of the IXth workshop ``Probability, Combinatorics and Geometry'' at Bellairs institute (2014) where this work started. In particular, we are grateful to Omer Angel and Simon Griffiths for interesting discussions. We also thank Frédéric Paulin for a question raised in 2008 which eventually yields to this work. Last we thank the referee for a relevant question which yields to Proposition \ref{remreferee}.

\bigskip

\begin{center} \hrulefill \textit{ In this paper, unless mentioned, we only consider bounded sequences $(a_{i})_{i\geq 1}$.} \hrulefill  \end{center}

\section{Tracking a uniform point}

\label{sec:uniformpoint}
The goal of this section is to give a necessary and sufficient condition for the height of a typical point of $\mathcal T_n$ (i.e. sampled according to the normalized  length measure $\mu_n$) to converge in distribution towards a finite random variable. For bounded sequence $(a_{i})_{i \geq 1}$ this condition is just
$$
\mathsf{H}( \mathbf{a}) = \sum_{i=1}^{+\infty} \frac{a_i^2}{A_i} < \infty.
$$
We will more precisely show that the above display is a necessary and sufficient condition for the convergence of the random measure $\mu_{n}$ towards a random probability measure $\mu$ carried by the limiting tree $ \mathcal{T}$. We begin by introducing a piece of notation.

\subsection{Notation}
\label{sec:not}

\noindent \textbf{$\mathbb{R}$-trees as subsets of $  \ell^{1}( \mathbb{R})$.} We briefly recall here some definitions about  $\mathbb R$-trees and refer to \cite{Eva08,LG06} for precisions. An $\mathbb{R}$-tree is a metric space $(\mathcal T, \delta)$ such that for every $
x,y\in \mathcal T$, there is a unique arc from $x$ to $y$ and this arc is isometric to a segment in $ \mathbb{R}$.  If $a,b \in \mathcal{T}$ we denote by $ \llbracket a,b\rrbracket$ the geodesic line segment between $a$ and $b$ in $ \mathcal{T}$. 
The degree (or multiplicity) of a point $x \in \mathcal{T}$ is the number of connected components of $ \mathcal{T}\backslash \{x\}$. A point of degree $1$ is a called a leaf and a point of degree at least $3$ is called a branch point. \medskip

Let $ \mathbf{a} = (a_{i})_{i \geq 1}$ be a sequence of positive reals, and $A_{i} = a_{1}+ ... + a_{i}$, for $i \geq 1$, the associated sequence of partial sums. From $\mathbf a$, we build a sequence of random trees  $(\mathcal T_n)_{n \geq 1}$ by grafting randomly closed segments (also called branches) of lengths $a_{i}, i \geq 1$ inductively as described in the introduction. To be more precise, we follow the initial approach of Aldous \cite{Ald91a} and build $ \mathcal{T}_{n}$ as a subset of $ \ell^1( \mathbb{R})$. The tree $ \mathcal{T}_{1}$ is $\{ (x,0,0,\ldots) : x \in [0,a_{1}]\}$ and recursively for every $n \geq 1$, conditionally on $ \mathcal{T}_{n}$, we pick $ (u^{(n)}_{1}, \ldots , u_{n}^{(n)}, 0, 0, \ldots) \in \mathcal{T}_{n}$ a uniform point on $ \mathcal{T}_{n}$ and set 
$$ \mathcal{T}_{n+1} := \mathcal{T}_{n} \cup \big\{ (u^{(n)}_{1}, \ldots , u_{n}^{(n)}, x, 0 , 0, \ldots) \in \ell^1( \mathbb{R}) : x \in [0, a_{n+1}]\big\}.$$
The point $\rho = (0,0, \ldots)$ will be seen as the root of the trees $ \mathcal{T}_{n}$. With this point of view, the trees $ \mathcal{T}_{n}$ are increasing closed subsets of $ \ell^1 ( \mathbb{R})$ and we can define their increasing union  $$\mathcal T^*=\bigcup_{n\geq 1} \mathcal T_n.$$ 
Note that $ \mathcal{T}^* \subset \ell^1( \mathbb{R})$ will not be closed in general (or equivalently complete). We let  $ \mathcal{T}$ denote its closure (or completion), which is therefore a random closed subset of $ \ell^1( \mathbb{R})$. For us, $ \mathcal{T}$ and $ \mathcal{T}_{n}$ once endowed with their length metric $ \delta$, will be viewed as random $ \mathbb{R}$-trees (recall that, in general, the completion of an $\mathbb R$-tree is an $\mathbb R$-tree -- see e.g. \cite{Imr77}). In the rest of this article, we will be loose on the fact that $ \mathcal{T}_{n}, \mathcal{T}$ are subsets of $ \ell^1( \mathbb{R})$ and will use it only when necessary for technical proofs.

\bigskip

\noindent \textbf{General notation.} Let  $ (\mathcal{F}_{n})_{n \geq 1}$ denote the associated filtration generated by $ ( \mathcal{T}_{n})_{n \geq 1}$,  and write $ \mathsf{b}_{i}$ for the segment or branch of index $i$ which is seen as a subset of $ \mathcal{T}_{n}$ for each $n \geq i$.  A moment of thought shows that $ \mathcal{T} \backslash \mathcal{T}^*$ is only made of leaves of $ \mathcal{T}$.  
We should stress that, although our main goal is to study some geometric properties of the sole tree $\mathcal T$, we will often need to work with its subtrees $\mathcal T_n,n\geq 1$. In that aim, we label the leaves of $\mathcal T^*$ by order of apparition in the aggregation procedure, so that when observing $\mathcal T$, we also know $\mathcal T_n$, which is simply the subtree of $\mathcal T$ spanned by the root and the leaves labeled $1,\ldots,n$, $\forall n \geq 1$. This property is automatic when $ \mathcal{T}_{n}$ is constructed as a subset of $ \ell^{1}( \mathbb{R})$ as before since the $i$th branch ranges over the $i$th coordinate of $ \ell^{1}( \mathbb{R})$.

Besides, as already  mentioned, we denote by ${\mu}_{n}$ the length measure on $ \mathcal{T}_{n}$  normalized by $A_{n}^{-1}$ to make it a probability measure. Also, to lighten notation, we write $ \mathsf{ht}(x) = \delta (x, \rho)$ for the height of  $x \in \mathcal{T}$.  

Thanks to the nested structure of the trees $(\mathcal T_n)_{n\geq 1}$, for $k \geq 1$ and for any point $ x \in \mathcal{T}$, we can make sense of $[x]_{k}$ the projection of $x$ onto $ \mathcal{T}_{k}$, that is the (unique) point of $ \mathcal{T}_{k}$ that minimizes the distance to $x$. If $A \subset \mathcal{T}$, for all $n\geq i$ we denote by 
 \begin{equation} \label{eq:subtrees}
 \mathcal{T}_{n}^{(i)}(A) = \big\{ x \in \mathcal{T}_{n} : [x]_{i} \in A\big\},
  \end{equation} 
 the subtree ``descending from'' $A$ in $ \mathcal{T}_{n}$.
Similarly we let $ \mathcal{T}^{(i)}(A) = \{ x \in \mathcal{T} : [x]_{i} \in A\}$, the subtree ``descending from'' $A$ in $ \mathcal{T}$. Note that these definitions depend in general on the integer $i$. E.g., $$\mathcal{T}^{(2)}(\mathcal T_1)  \subsetneq \mathcal{T} ^{(1)}(\mathcal T_1)=\mathcal T.$$

\medskip

\noindent \textbf{Stems.} A stem of a tree is a maximal open segment that contains no branch point. We will use a genealogical labeling of the stems of the trees $ ( \mathcal{T}_{n})_{n \geq 1}$ by the ternary tree $$ \mathcal{G} =  \bigcup_{i \geq 0} \{0,1,2\}^i,$$ with the usual genealogical order $\preccurlyeq$.  Formally the first branch $ \mathsf{b}_{1}$ is labeled by $\varnothing$. Once we graft a branch on it, it is split into three stems denoted (arbitrary) by $0,1,2$. Recursively, when the stem labeled $u \in \mathcal{G}$ is split into three by grafting a new branch on it, we denote $u0, u1,u2$ the three stems created. Here and later we implicitly identify a stem with its label.  When $ \mathcal{T}_{n}$ is built after $n$ graftings we denote by $ \mathcal{G}_{n} \subset \mathcal{G}$ the set of all stems of $ \mathcal{T}_{n}$.

When $u \in \mathcal{G}_{i}$ is a stem of $ \mathcal{T}_{i}$ we lighten the notation introduced in \eqref{eq:subtrees} and set 
$$ \mathcal{T}_{n}(u) := \mathcal{T}_{n}^{(i)}(u) \quad \mbox{ and } \quad  \mathcal{T}(u) := \mathcal{T}^{(i)}(u).$$ It is easy to check that these definition do not depend on $i$ when $u$ happens to belong to several $ \mathcal{G}_{i}$. The last remark is also valid if $u$ is the closure of a stem. We use the notation $  \mathsf{L}({u})$ for the length of the stem $u$ and introduce for $u \in \mathcal{G}_{n}$ 
$$  \mathbf{a}(u) = (a_{i}(u))_{ i \geq 1} = (0)_{1 \leq i \leq n-1} \cup \{\mathsf{L}({u})\} \cup \left(a_{i} \mathbbm{1}_{\left\{a_{i}  \mathrm{\  is\  grafted\ on\  } \mathcal{T}_{i-1}(u)\right\}}\right)_{i \geq n+1}$$ for the sequence of lengths of branches that are recursively grafted onto the stem $u$ or its descendants, with the convention that the first branch is the stem $u$ appearing at time $n$. Note that $\mathbf a(u)$ corresponds to the lengths of branches used to construct $\mathcal T(u)$. We will sometimes need to consider a notion of height in these subtrees.  Let $\overline u=u\cup\{a_u\}\cup\{b_u\}$ be the closure of $u$ in $\mathcal T$, where $a_u$ designs the vertex closest to the root. Then we define the height of a vertex $x \in \mathcal T(u)$ as the distance $\delta(a_{u},x)$ and the height of the tree $\mathcal T(u)$ as  the supremum of the distances  $\delta(a_{u},x)$
when $x$ runs over $\mathcal T(u)$.

\begin{remark} Almost surely the set of branch-points of $ \mathcal{T}$ is dense in $\mathcal T$. Indeed, since the sequence $(a_i)_{i \geq 1}$ is bounded, $A_i \leq ci$ for some constant $c<\infty$ and all $i$. In particular   $$\sum_{i \geq 1}  \frac{1}{A_{i}}=\infty$$ and the  Borel--Cantelli lemma implies that  infinitely many branches will be grafted on each stem, almost surely. If $a_{i} \to 0$ we even have that the set of leaves of $ \mathcal{T}$ is dense in $ \mathcal{T}$ a.s.. \end{remark}

\subsection{Height of a random point} 
\label{sec:height}
We begin with a simple key observation. Let $n \geq 2$ and conditionally on $ \mathcal{T}_{n}$ pick a point $Y_{n}$ uniformly distributed according to the measure $\mu_{n}$. Two cases may happen:
\begin{itemize} 
\item with probability $1 - a_{n}/A_{n}$: the point $Y_{n}$ belongs to the tree $ \mathcal{T}_{n-1}$, that is $[Y_{n}]_{n-1}= Y_{n},$ and conditionally on this event $[Y_{n}]_{n-1}$ is uniformly distributed over $ \mathcal{T}_{n-1}$,
\item with probability $a_{n}/A_{n}$: the point $Y_{n}$ is located on the last branch $\mathsf b_n$ grafted on $ \mathcal{T}_{n-1}$. Conditionally on this event, $Y_{n}$ is uniformly distributed on this branch and its projection $ [Y_{n}]_{n-1}$ on the tree $ \mathcal{T}_{n-1}$ is independent of its location on the $n$th branch and is uniformly distributed on $ \mathcal{T}_{n-1}$, given $\mathcal T_{n-1}$.
\end{itemize}
From this observation we deduce that  $( \mathcal{T}_{n-1},[Y_{n}]_{n-1}) = ( \mathcal{T}_{n-1}, Y_{n-1})$ in distribution and more generally, $( \mathcal{T}_{k}, [Y_{n}]_{k}) = ( \mathcal{T}_{k}, Y_{k})$ in distribution for all $ 1 \leq k \leq n$. Note however an important subtlety: given the tree $\mathcal T_n$, the point $[Y_{n}]_{n-1}$ \textit{is not} uniformly distributed on its subtree $\mathcal T_{n-1}$ since $[Y_{n}]_{n-1}$ is located on a branch point of $ \mathcal{T}_{n}$ with probability $a_{n}/A_{n}$.

Reversing the process, it is possible to build a sequence $(\mathcal{T}_{n},X_{n})_{n \geq 1}$ recursively such that $[X_{n}]_{k} = X_{k}$ for all $k \leq n$ and such that $ ( \mathcal{T}_{n},X_{n}) = ( \mathcal{T}_{n}, Y_{n})$ in law for every $n$. 
To do so, consider an independent sample $(U_{i},V_i,i\geq 1)$ of i.i.d.~uniform random variables on $(0,1)$. Let first $\mathcal T_1$ be a segment of length $a_1$, rooted at one end, and let $X_1$ be the point on this segment at distance $a_1 V_1$ from the root.  We then proceed recursively and assume that the pair $(\mathcal T_n,X_n)$ has been constructed. Then: 
\begin{enumerate}
\item[$\bullet$] if $U_{n+1} \leq a_{n+1}/A_{n+1}$, we branch a segment of length $a_{n+1}$ on $X_n$ to get $\mathcal T_{n+1}$ and let $X_{n+1}$ be the point on this segment at distance $a_{n+1}V_{n+1}$ from the branchpoint $X_n$, 
\item[$\bullet$] if $U_{n+1} > a_{n+1}/A_{n+1}$, we branch a segment of length $a_{n+1}$ at a point chosen uniformly (and independently of $X_{n}$) at random in $\mathcal T_n$, and set $X_{n+1}=X_n$. 
\end{enumerate}
Clearly, $[X_{n}]_{k} = X_{k}$ for $1 \leq k \leq n$ and it is easy to see by induction that $ ( \mathcal{T}_{n},X_{n})$ and $( \mathcal{T}_{n}, Y_{n})$ have the same distribution for all $n\geq 1$. It is important to notice that in this coupling, the distance between $X_n$ and the root $\rho$ is non-decreasing, and more precisely that for any $n \geq m \geq 0$,
  \begin{eqnarray}  \label{eq:distance}  
\delta( X_{n}, \mathcal{T}_{m}) = \delta(X_{n}, X_{m}) = \sum_{i= m+1}^n a_i V_i \mathbbm 1_{\left\{U_i \leq \frac{a_i}{A_i}\right\}},
 \end{eqnarray}
where we have set $X_{0}= \mathcal T_0= \rho$. Recalling the definition of $ \mathsf{H}( \mathbf{a})$ we see that $ \lim_{n \to \infty} \mathbb{E}[ \mathsf{ht}(X_{n})] = \mathsf{H} (\mathbf{a})/2$. Therefore, when $ \mathsf{H}( \mathbf{a}) < \infty$,  the sequence ($  \mathsf{ht}(X_{n})$) converges  and moreover $(X_{n})$ is a Cauchy sequence, by (\ref{eq:distance}), almost surely. So, in this case, $(X_{n})$ converges a.s. in $\mathcal T$, by completeness.  The converse is also true:

\begin{proposition}[Finiteness of a typical height] \label{prop:unif} For bounded sequences $(a_{i})_{i \geq 1}$,
$$ (X_{n}) \text{ converges in } \ \mathcal T \text{ a.s.}\quad \iff \quad \mathsf{H}( \mathbf{a}) < \infty.$$
Moreover, when $\mathsf{H}( \mathbf{a}) < \infty$, if $X:=\lim_{n\rightarrow \infty} X_n$, we have 
 $$  \mathbb{E}\big[e^{\lambda \mathsf{ht}(X)}\big]  \leq e^{ \lambda \mathsf{H}( \mathbf{a})}, \quad  \text{ for all } \lambda \in \Big[0,\big(\textstyle{\sup_{i\geq 1} a_{i}}\big)^{-1}\Big]. $$
 \end{proposition}

 \proof By \eqref{eq:distance}, the convergence of $(X_{n})$ is equivalent to the convergence of the series $ \sum_i a_{i}V_{i} \mathbbm{1}_{\{U_{i} \leq a_{i}/A_{i}\}}$ and so the first point follows from the classical  three series theorem.  To establish the exponential bound, note that for all $n \geq 1$,
  \begin{eqnarray*} \mathbb{E}\Big[e^{\lambda  \mathsf{ht}(X_{n})}\Big] &=& \prod_{i=1}^n \left( \frac{A_{i}-a_{i}}{A_{i}} + \frac{a_{i}}{A_{i}}  \mathbb{E}\left[e^{\lambda a_{i}  V_{i}}\right] \right)\\
  &=& \prod_{i=1}^n \left( \frac{A_{i}-a_{i}}{A_{i}} + \frac{a_{i}}{A_{i}}  \frac{1}{\lambda a_{i}}\left( e^{\lambda a_{i}}-1\right) \right). \end{eqnarray*}
 Then, since $\lambda a_{i} \leq 1$, we can use the bound $e^x \leq 1 +x +x^2$  valid for all $x \in [0,1]$, and also $\log(1+x) \leq x$ for $x \geq 0$, to get
 \begin{eqnarray*} \prod_{i=1}^n \left( \frac{A_{i}-a_{i}}{A_{i}} + \frac{a_{i}}{A_{i}}  \frac{1}{\lambda a_{i}}\left( e^{\lambda a_{i}}-1\right) \right) & {\leq} & \prod_{i=1}^n \left( \frac{A_{i}-a_{i}}{A_{i}} + \frac{a_{i}}{A_{i}}  (1 + \lambda a_{i}) \right) \\ & =& \exp\left( \sum_{i=1}^n \log\left(1 + \lambda \frac{a_{i}^2}{A_{i}}\right) \right) \\
  & {\leq} & \exp\left( \lambda \sum_{i=1}^n \frac{a_{i}^2}{A_{i}}\right).
   \end{eqnarray*}
   Letting $n \to \infty$ we get the desired bound.
 \endproof

\begin{remark} \label{rem:leaf} By equation (\ref{eq:distance}) we get that $ \mathbb{P}(X_{n}=X_{n_{0}}, \forall n \geq n_{0}) = A_{n_{0}}/A_{\infty}$ and so, with probability one, the sequence $(X_n)$ is eventually constant if and only if $\sum_i a_i$ is convergent.
\end{remark}
 \begin{remark} In the case of unbounded sequences $(a_{i})_{i \geq 1}$ (not considered in this paper) the three series theorem shows that $(X_{n})$ converges a.s.~iff there exists some $ \varepsilon>0$ such that 
 $$  \sum_{i\geq 1} \frac{a_{i}}{A_{i}} \mathbbm{1}_{\{a_{i} \geq \varepsilon\}} < \infty \quad \mbox{ and } \quad \sum_{i \geq 1} \frac{a_{i}^2}{A_{i}} \mathbbm{1}_{\{ a_{i} \leq \varepsilon\}} < \infty.$$
 \end{remark}

\medskip

\noindent \textbf{Examples:}

\begin{enumerate}
\item If the sum $\sum_{i\geq 1} a_i$ is finite, or if $a_{i} \leq i^{-\eps + \circ(1)}$ for some $\eps>0$, then $  \mathsf{H}( \mathbf{a})$ is finite (see Lemma \ref{lem:manip} (ii)) and so the tree $ \mathcal{T}_{n}$ has a typical height which remains bounded as $ n \to \infty$. Proposition \ref{prop:compact} and Proposition \ref{prop:finitelength} actually state that in these cases the maximal height  of the tree $\mathcal{T}_{n}$ remains bounded as $ n \to \infty$. 
\item If  $a_i \sim (\ln i)^{-\lambda}$ for some $\lambda \leq 1$ then $ \mathsf{H}( \mathbf{a}) = \infty$ and so the typical height of $ \mathcal{T}_{n}$ blows up. On the other hand, if $a_i \sim (\ln i)^{-\lambda}$ for some $\lambda >1$ then $ \mathsf{H}( \mathbf{a}) < \infty$ and the typical height of $ \mathcal{T}_{n}$ thus remains bounded. In this case, we do not know whether the maximal height of $ \mathcal{T}_{n}$ remains stochastically bounded as $n \to \infty$.

\item Consider the sequence 
$$a_i = i^{-1/2} + 
\mathbbm 1_{\{i \in \mathbb N^3\}} \qquad \forall i \geq 1.$$
Clearly, $A_i \sim 2 \sqrt i$ \ and $\mathsf{H}( \mathbf{a})<\infty$. Although the typical height of $ \mathcal{T}_{n}$ remains bounded, the tree $\mathcal T$ is not compact since it contains an infinite number of  branches of length greater than $1$. (In fact, this tree is even unbounded, see Subsection \ref{sec:bounded}.)
\end{enumerate}

\subsection{Convergence of  the length measure $\mu_n$}
\label{section:cvmun}

By Proposition \ref{prop:unif}, when $ \mathsf{H}( \mathbf{a}) = \infty$  the height of a random point in $ \mathcal{T}_{n}$ sampled according to $\mu_{n}$ tends in probability to $\infty$. It follows that the sequence of probability measures $(\mu_{n})$ cannot converge weakly in this context. However we will see that it does converge as soon as $ \mathsf{H}(\mathbf a) < \infty$. With no loss of generality, we assume in the sequel that the tree $\mathcal T$ is built jointly with the sequence $(X_n)$, as explained in the previous section.

\begin{theorem}[Convergence of the length measures]  \label{thm:measures}Suppose that $ \mathsf{H}( \mathbf{a}) < \infty$. Then almost surely,  there exists a probability measure $\mu$ on $ \mathcal{T}$ such that 
$$\mu_{n} \to \mu \quad  \mbox{weakly as } n \to \infty.$$
Furthermore, conditionally on $\mu$, the point $X=\lim_{n \rightarrow \infty} X_n$ is distributed according to $\mu$ almost surely and there is the dichotomy:
\begin{itemize}
\item if $\sum_i a_{i} = \infty$ then $\mu$ is a.s.\,\,supported by the leaves of $ \mathcal{T}$,
\item if $\sum_i a_{i} < \infty$ then $\mu$ is a.s.\,\,supported by the skeleton of $ \mathcal{T}$ and coincides with the normalized length measure of $ \mathcal{T}$. 
\end{itemize}
\end{theorem}

\medskip

To get a precise meaning of this theorem, recall that the trees $ \mathcal{T}_{n}, n \geq 1$ and $ \mathcal{T}$ were actually constructed as closed subsets of $  \ell^{1}( \mathbb{R})$. Hence, the random probability measures $ \mu_{n}$ are just random variables with values in the Polish space of probability measures on $\ell^1(  \mathbb{R})$ endowed with the Lévy-Prokhorov distance (which induces the weak convergence topology). Recall that the Lévy-Prokhorov distance on the probability measures of a metric space $(E,d)$ is given by 
$$ \mathrm{d_{LP}}(\mu, \nu) = \inf\Big\{ \varepsilon>0 : \nu(A) \leq \mu(A^{( \varepsilon)})+ \varepsilon \mbox{ and }\mu(A) \leq \nu(A^{( \varepsilon)})+ \varepsilon, \ \mbox{ for all Borel }A \subset E \Big\},$$ and where $A^{( \varepsilon)} = \{ y \in E : d(y,A) \leq \varepsilon\}$ is the $ \varepsilon$-enlargement of $A$. 

\medskip

\begin{proposition}
\label{remreferee}
Let $\mu_{n,\mathrm{leaves}}$ be the empirical measure on the $n$ leaves of $\mathcal T_n$. When $\mathsf H(\mathbf a)<\infty$
$$
\mu_{n,\mathrm{leaves}}\to \mu \quad  \mbox{weakly as } n \to \infty, \quad a.s.
$$
where $\mu$ is the probability measure arising  in Theorem \ref{thm:measures}.
A similar result holds for the empirical measures on the branch points of $\mathcal T_n$, or on the set of leaves and branch points of $\mathcal T_n$.
For the Brownian CRT, this implies that the measure $\mu$ corresponds to the usual uniform measure carried by this tree.
\end{proposition}

\medskip

The proofs of Theorem \ref{thm:measures} and Proposition \ref{remreferee} occupy the rest of this subsection. To prove the first point of the theorem we will show that $(\mu_{n})$ is a Cauchy sequence.  We point out that this is not a direct consequence of Proposition \ref{prop:unif}. Indeed, as noticed in the previous section, given the tree $\mathcal T$, the variable $X_n$ is \emph{not} distributed according to $\mu_n$ since it is equal to a branch point of $\mathcal T$ with a strictly positive probability. 
 We start by introducing a family of martingales which will play an important role.

\medskip

\noindent \textbf{Mass martingales.} Let $ C \subset \mathcal{T}_{i}$ be measurable for  $ \mathcal{F}_{i}$ and  recall the notation $ \mathcal{T}^{(i)}_{n}(C)$ for $n \geq i$ and $  \mathcal{T}^{(i)}(C)$ introduced in (\ref{eq:subtrees}).  Set then $M_{n}(C) = \mu_{n}( \mathcal{T}^{(i)}_{n}(C))$ to simplify notation. Since the branches are grafted uniformly on the structure at each step, we have conditionally on $ \mathcal{F}_{n}$
$$ \left\{ \begin{array}{lcc} \vspace{0.2cm}\displaystyle M_{n+1}(C) = (A_{n} \cdot  M_{n}(C) + a_{n+1})/A_{n+1} & \mbox{ with proba.} & M_{n}(C),\\
\displaystyle M_{n+1}(C) = A_{n} \cdot  M_{n}(C)/A_{n+1}& \mbox{ with proba.} & 1- M_{n}(C).\end{array}\right.$$ It readily follows that $(M_{n}(C))_{n\geq i}$ is a martingale with respect to $ (\mathcal{F}_{n})_{n \geq i}$ and since it takes values in $[0,1]$, it converges almost surely to its limit $M(C) \in [0,1]$. This limit  $M(C)$ is the natural candidate for the value of $\mu( \mathcal T^{(i)}(C))$ of the possible limit $\mu$ of  $(\mu_n)$.

\begin{remark}[Generalized Polya urn] 
\label{remark:Polya}
These martingales are also known as ``generalized Polya urns'' in the theory of reinforced processes. In general, it is a subtle question to discuss whether $M(C)$ can have atoms in $\{0,1\}$, see \cite{Pem90}. However, in our context, since the sequence $(a_i)_{i \geq 1}$ is bounded, it follows from Pemantle's work \cite{Pem90} that $M(C) \in (0,1)$ almost surely when $C$ and $C^c$ have positive length measures. Let us emphasize an important consequence for us. Consider $C \subset \mathcal{T}_{i}$ with positive length measure and $ \mathcal{F}_{i}$-measurable and let $J$ be an infinite subset of $\mathbb N$. Then, 
$$
\sum_{j \in J, j \geq i} M_j(C)=\infty \quad \text{a.s.}
$$
and the conditional version of the Borel--Cantelli lemma implies that almost surely an infinite number of branches $\mathsf b_j, j \in J$ belong to the subtree $\mathcal T^{(i)}(C)$.
\end{remark}

\begin{lemma}  \label{lem:tight} Assume  $\mathsf{H}( \mathbf{a}) < \infty$. Then almost surely, for any $ \varepsilon>0$, there exists (a random) $n_{0}$ such that 
$$ \mu_{n}\big( \mathcal{T}_{n_{0}}^{ (\varepsilon)}\big) \geq 1- \varepsilon \quad \text{ for all } n\geq 1.$$
\end{lemma}
\proof 
We use the construction of $( \mathcal{T}_{n},X_n)$ of Section \ref{sec:height}. Fix $\eps>0$ and a (deterministic) integer $n_0$ and consider the stopping time (with respect to the filtration $(\mathcal F_n)$) defined by
 $$ \theta = \inf \left\{ n \geq 1 : \mu_{n}( \mathcal{T}_{n_{0}}^{( \varepsilon)}) < 1- \varepsilon \right\}.$$
Note that
\begin{eqnarray*}
 \mathbb P \left(\theta<\infty, \delta (X_{\theta}, X_{n_0}) \leq \eps\right) &=& \sum_{n \geq 1} \mathbb E\big[ \mathbb P\left(\theta=n, \delta (X_{n}, X_{n_0}) \leq \eps | \mathcal F_n\right) \big] \\
 &\leq & \sum_{n \geq 1} \mathbb E\left[ \mathbbm 1_{\{\theta=n\}} \right](1-\eps) = (1-\eps) \mathbb P(\theta<\infty)
\end{eqnarray*}
where we have used that the distribution of $X_n$ given $\mathcal F_n$ is $\mu_n$, as well as the definition of $\theta$, to get the second inequality. This yields
 \begin{eqnarray*}  \varepsilon \cdot \mathbb{P}( \theta < \infty)& \leq & \mathbb P\left(\theta<\infty, \delta (X_{\theta},X_{n_0}) >\eps \right)\\ 
 & \leq&  \mathbb P\left( \delta (X,X_{n_0}) >\eps \right) \\
 & \underset{ \eqref{eq:distance}}{\leq} & \frac{1}{2\eps}\sum_{i=n_{0}+1}^\infty \frac{a_{i}^2}{A_{i}}.  
 \end{eqnarray*} 
Since the right-hand side can be made arbitrarily small by letting $n_{0} \to \infty$, we get that almost surely, for every $ \varepsilon>0$ (rational say),  there exists (a random) $n_{0} \geq 1$ such that 
$ \mu_{n}\big( \mathcal{T}_{n_{0}}^{ (\varepsilon)}\big) \geq 1- \varepsilon$ for all $n \geq 1$.
\endproof

\begin{lemma} 
\label{lem:Cauchy}
Assume  $\mathsf{H}( \mathbf{a}) < \infty$.  Then almost surely $(\mu_{n})$ is a Cauchy sequence for the Lévy-Prokhorov distance.
\end{lemma}
\proof 
For any $0 \leq k \leq n$, let $[\mu_{n}]_{k}$ be the measure $\mu_{n}$ projected onto $ \mathcal{T}_{k}$, that is the push forward of $ \mu_{n}$ by $ x \mapsto [x]_{k}$. The following assertions hold almost surely. Fix $ \varepsilon>0$, it follows from the last lemma that there exists (a random) $n_{0}$ such that for all $n \geq 1$
 \begin{eqnarray} 
 \mathrm{d_{LP}}( \mu_{n} ; [\mu_{n}]_{n_{0}}) \leq \varepsilon.   \label{eq:etape1}\end{eqnarray}
Indeed, if $Y_{n}$ is sampled according to $\mu_{n}$ then we have $\delta( Y_{n}, [Y_{n}]_{n_{0}}) \leq \varepsilon$ with probability at least $1- \varepsilon$. Since $[Y_{n}]_{n_{0}}$ is distributed as $[\mu_{n}]_{n_{0}}$ this readily implies the (\ref{eq:etape1}). We then decompose $ \mathcal{T}_{n_{0}}$ into a finite number of $ \mathcal{F}_{n_{0}}$-measurable pieces $C_{1}, \ldots , C_{K}$ of diameter less than $ \varepsilon$ (note that $K$ is random). For each of these pieces recall the definition of the martingale $M_{n}(C_{j})$ for $n \geq n_{0}$. In particular with our notation we have $M_{n}(C_{j}) = [\mu_{n}]_{n_{0}}(C_{j})$. Next, note that when 
$$\sum_{i =1}^K |M_{n}(C_{i})- M_{m}(C_{i})| \leq   \varepsilon,$$ we can couple $X \sim [\mu_{n}]_{n_{0}}$ and $X' \sim [\mu_{m}]_{n_{0}}$ so that $X$ and $X'$ belong to the same set $C_{i}$ with probability at least $1- \varepsilon$. This implies that $ \mathrm{d_{LP}}([\mu_{n}]_{n_{0}} ; [\mu_{m}]_{n_{0}}) \leq \varepsilon$. Since the martingales $(M_{n}(C_{i}))$ converge as $n \to \infty$, the last display is eventually fulfilled for $n,m$ large enough. As a result, for $n,m$ large enough
$$\mathrm{d_{LP}}([\mu_{n}]_{n_{0}} ; [\mu_{m}]_{n_{0}}) \leq \varepsilon.$$
Combining the last display with \eqref{eq:etape1} we get that for $n,m$ large enough, $ \mathrm{d_{LP}}( \mu_{n} ; \mu_{m}) \leq 3 \varepsilon$. Hence $(\mu_{n})$ is almost surely Cauchy for the Lévy--Prokhorov distance on $ \ell^{1}( \mathbb{R})$.\endproof 

\medskip

\proof[Proof of Theorem \ref{thm:measures}] The existence of the almost sure limit $\mu$ of $(\mu_n)$ is ensured by the previous lemma.

\smallskip

\noindent \textsc{Distribution of $X$.} Recall from Section \ref{sec:height} that $X_{n} \sim \mu_{n}$, {given $\mu_n$}. In particular, for any $n \geq m$,  $X_{m} = [X_n]_{m}$ is distributed according to  $[\mu_{n}]_{m}$, {given $\mu_n$}. Letting $n \to \infty$ {and using the continuity of the projection on $\mathcal T_m$ for the Lévy-Prokhorov distance}, we obtain that $X_{m} \sim [\mu]_{m}$, {given $\mu$}. Now, let $m\rightarrow \infty$. On the one hand, {according to the arguments developed in the proof of Lemma \ref{lem:Cauchy}}, $[\mu]_{m} \to \mu$ almost surely for the Lévy-Prokhorov metric. On the other hand, $X_{m} \to X$ almost surely. It follows that $X \sim \mu$ {almost surely given $\mu$}.

\noindent \textsc{Support of $\mu$.}  Since  $X \sim \mu$ almost surely given $\mu$, we only need to show that $ \mathbb{P}( X \mbox { is a leaf  of } \mathcal{T}) =1$ or $0$ according to $\sum_i a_{i} = \infty$ or $\sum_i a_{i} < \infty$. By the construction of $X_{n}$ and $X$, we have 
$$ \mathbb{P}( X \mbox{ is a leaf in } \mathcal{T}) = \lim_{m \to \infty} \lim_{n \to \infty}\mathbb{P}(X_{n} \notin \mathcal{T}_{m}).$$
If $ \sum_i {a}_{i} = \infty$, by Remark \ref{rem:leaf}, the sequence $(X_{n})$ escapes from any finite tree $ \mathcal{T}_{m}$ almost surely   and so $ \mathbb{P}( X \mbox{ is a leaf in } \mathcal{T}) =1$. Conversely if $ \sum_i a_{i} < \infty$, then $X_{n}=X$ eventually  so $ \mathbb{P}( X \mbox{ is a leaf in } \mathcal{T}) ~=~0$. In this case, $(\mu_{n})$ converges clearly towards the normalized length measure on $\mathcal{T}$.
\endproof 

\medskip

\proof[Proof of Proposition \ref{remreferee}] For all $i$ and then all $ \mathcal{F}_{i}$-measurable  $ C\subset \mathcal{T}_{i}$, we let $L_n(C)$ be the number of leaves in $ \mathcal{T}^{(i)}_{n}(C)$,  $n>i$. Note that  $\mu_{n,\mathrm{leaves}}(\mathcal{T}^{(i)}(C))=L_n(C)/n.$ Omitting easy details, we claim that the almost sure convergence $\mu_{n,\mathrm{leaves}} \rightarrow \mu$ will be proved if we check that  $n^{-1} L_n(C) \rightarrow \mu(\mathcal{T}^{(i)}(C))$ a.s., for all $ \mathcal{F}_{i}$-measurable $ C \subset \mathcal{T}_{i}$, for all $i$.  So fix such a couple $(C,i)$ and observe that
$$
N_n(C):=L_n(C)-\sum_{k=i}^{n-1}M_k(C), \quad n >i
$$
defines a centered martingale, such that $|N_{n+1}(C)-N_n(C)|\leq 1$, a.s for all $n>i$. Applying Azuma-Hoeffding inequality, we get that for all $\varepsilon>0$
$$
\mathbb P\left(|N_n(C) |\geq \eta n^{\frac{1}{2}+\varepsilon}\right) \leq 2 \exp\left(- \frac{\eta^2n^{1+2\varepsilon}}{2(n-i)}\right), \quad \forall n >i.
$$
By Borel-Cantelli's lemma, this obviously implies that 
$$
\frac{N_n(C)}{n^{\frac{1}{2}+\varepsilon}} \underset{n \rightarrow \infty}{\overset{\mathrm{a.s.}}\longrightarrow} 0
$$
for all $\varepsilon>0$ and in particular that $n^{-1} N_n(C)\rightarrow 0$ a.s. On the other hand, when $\mathsf H(\mathbf a)<\infty$,  Theorem \ref{thm:measures} implies that $M_n(C) \rightarrow \mu(\mathcal{T}^{(i)}(C))$ a.s., and so the Ces\`aro mean $n^{-1}\sum_{k=i}^{n-1}M_k(C) \rightarrow \mu(\mathcal{T}^{(i)}(C))$ a.s. as well. Hence $n^{-1} L_n(C)\rightarrow \mu(\mathcal{T}^{(i)}(C))$ a.s. as expected. The proof holds similarly for the empirical measure on the branch points.
\endproof

\subsection{Boundedness of the whole tree}
\label{sec:bounded}
By Proposition \ref{prop:unif}, if the tree $\mathcal T$ is bounded we must have $\mathsf{H}( \mathbf{a}) < \infty$. We refine this a little: 
 
 \begin{proposition} A necessary condition for the tree $ \mathcal{T}$ to be bounded is that $a_{i} \rightarrow 0$ as $i \rightarrow \infty$.
 \end{proposition}

\proof To see this, assume that there is a real number $ \varepsilon>0$ and an infinite subset $J$ of $\mathbb N$ such that $a_i \geq  \varepsilon$ for all $i \in J$ (recall that  the $a_{i}$ are however supposed to be bounded). For each $i \in J$, let $\mathsf b_i^+$ denote the half part of the branch $\mathsf b_i$
composed by the points at distance at least $ \varepsilon/2$ from the vertex of $\mathsf b_i$ which is the closest to the root of $\mathcal T$.
Then, by an argument  similar  to that of Remark \ref{remark:Polya}, we know that almost surely, for each $\mathsf b_i^+$, $i \in J$, there is an infinite number of branches $\mathsf b_j, j \in J$ that belong to its descending subtree. Iterating the argument, we see that there is a path in $\mathcal T$ containing an infinite number of disjoint segments of lengths all greater than or equal to $ \varepsilon/2$. Hence $ \mathcal{T}$ is unbounded.  
\endproof
Using a variation of the above argument we even get 

\begin{proposition} Almost surely, 
$$ \mathcal{T} \mbox{ is compact} \iff \mathcal{T} \mbox{ is bounded}.$$
\end{proposition}
\proof The implication $ \Rightarrow$ is deterministically true. Notice that the events $\{ \mathcal{T} \mbox{ is not compact}\}$ and  \linebreak $\{ \mathcal{T} \mbox{ is not bounded}\}$ are contained in the tail $\sigma$-algebra generated by the gluings and so have probability $0$ or $1$. We suppose thus that $ \mathcal{T}$ is almost surely non-compact and will prove that it is almost surely non-bounded. We need a little notation. Fix $n \geq m$, the set $\mathcal{T}_{n} \backslash \mathcal{T}_{m}$ is a forest (a finite family of trees) whose highest tree is denoted by $ \tau(m,n)$ (we add its root to make it complete). It is easy to see that conditionally on $ \mathcal{F}_{m}$, the tree $\tau(m,n)$ is grafted on a uniform point of $ \mathcal{T}_{m}$. By monotonicity the limit $$ \xi = \lim_{m \to \infty} \lim_{n \to \infty}  \mathsf{ht}( \tau(m,n)) \in[0,\infty]$$ exists and is independent of $ \mathcal{F}_{m}$ for any $m \geq 0$. By the zero-one law $\xi$ is thus deterministic. Assume by contradiction that $\mathcal T$ is bounded a.s. Then $\xi<\infty$ and we must have $\xi >0$, otherwise $ \mathcal{T}$ would be pre-compact hence compact by completeness. Moreover, there exists then an integer $k$ such that 
  \begin{eqnarray} \label{bound1} \mathbb{P}\big( \mathsf{ht}( \mathcal{T}_{k}) \geq \mathsf{ht}( \mathcal{T})- \xi/4\big) \geq 1/2.  \end{eqnarray}
We denote by $C_{k}$ the $ \mathcal{F}_{k}$-measurable part 
$$ C_{k} = \{ x \in \mathcal{T}_{k} : \delta( \rho,x) \geq \mathsf{ht}( \mathcal{T}_{k}) - \xi/4\}.$$
Then for any $m \geq 1,$ consider the stopping time $\theta(m) = \inf\{ n \geq m : \mathsf{ht}(\tau{(m,n)})> \xi/2\}$ which is almost surely finite by definition of $\xi$. We put $\theta^0 = k$ and $\theta^r$ the $r$-fold composition $\theta \circ ... \circ \theta(k)$ to simplify notation. Recalling that for any $i \geq 0$, conditionally on $ \mathcal{F}_{\theta^i}$, the tree $\tau( \theta^i, \theta^{i+1})$ is grafted on a uniform point of $ \mathcal{T}_{\theta^i}$ we get
  \begin{eqnarray} \mathbb{P}\left( \bigcap_{i=0}^\infty \left\{\tau{(\theta^i, \theta^{i+1})} \mbox{ is not grafted on } \mathcal{T}^{(k)}_{\theta^i}( C_{k})\right\}\right) = \mathbb{E}\left[ \prod_{i=0}^\infty \left(1- \mu_{ \theta^i}(\mathcal{T}^{(k)}_{\theta^i}(C_{k}) )\right)\right].   \label{bound2}\end{eqnarray}
Remark \ref{remark:Polya} shows that  $\mu_{ n}(\mathcal{T}^{(k)}_{n}(C_{k}) )$ is a.s.~bounded away from $0$ uniformly in $n$ and so the last display is equal to $0$. This leads to a contradiction with \eqref{bound1} since grafting $\tau(\theta^i, \theta^{i+1})$ onto $\mathcal{T}^{(k)}_{\theta^i}(C_{k})$ gives a tree with height strictly greater than $\mathsf{ht}(\mathcal T_k)+\xi/4$.
\endproof

We will see in the forthcoming Proposition \ref{prop:compact} and Proposition \ref{prop:finitelength} that sufficient conditions for the compactness of $\mathcal T$  are either that $a_i\leq i^{-\alpha+\circ(1)}$ for some $\alpha \in (0,1]$ or that the series $\sum_i a_i$ is convergent. But we do not have a necessary and sufficient condition for boundedness or equivalently compactness of the tree, hence the following question :

\begin{open} Find a necessary and sufficient condition for $ \mathcal{T}$ to be bounded.  \label{open:bounded}
\end{open}

As mentioned in the Introduction, this problem was solved by Amini et al. \cite{ADGO14} for decreasing sequences $\mathbf a$: in these cases, with probability one, the tree $\mathcal T$ is bounded if and only if $\sum_{i\geq 1} i^{-1}a_i <\infty$. Note that in general this condition cannot be sufficient for boundedness: in the Example 3 of Section \ref{sec:height} the sum $\sum_{i\geq 1} i^{-1}a_i $ is finite, but the corresponding tree is unbounded since $a_i$ does not converge to 0.

\section{Infinite length case}
The goal of this section is to prove Theorem \ref{thm:main}. We will first prove (under more general conditions than those of Theorem \ref{thm:main}) that $ \mathcal{T}$ is compact using a covering argument which will also give the upper bound $ \mathrm{dim_{H}}( \mathcal{T}) \leq 1/\alpha$. The lower bound on the Hausdorff dimension then follows from a careful study of the random measure $\mu$ introduced in Theorem \ref{thm:measures} and, again, is valid under more general conditions than those of Theorem \ref{thm:main}.

\subsection{Compactness and upper bound}
The main result of this subsection is the following:
\begin{proposition} \label{prop:compact} Assume that $a_{i} \leq i^{-\alpha+\circ(1)}$ for some $\alpha \in (0,1]$. Then, almost surely, the random tree $ \mathcal{T}$ is compact and its Hausdorff dimension is at most $\alpha^{-1}$.
 \end{proposition}
 We point out that we more generally know that the tree $\mathcal T$ is compact, with a set of leaves of Hausdorff dimension less than $\alpha^{-1}$, as soon as $a_{i} \leq i^{-\alpha+\circ(1)}$ for some $\alpha>0$. This  
 follows from the previous result, together with the forthcoming Proposition \ref{prop:finitelength}. 
 That being said, we focus in the rest of this subsection on the proof of Proposition \ref{prop:compact} and assume that $a_{i} \leq i^{-\alpha+\circ(1)}$ for $\alpha \in (0,1]$. We note with Lemma \ref{lem:manip} (ii) that this implies that
$$
\sum_{i=n}^{\infty} \frac{a_i^2}{A_i} \ \leq \ n^{-\alpha+\circ(1)},
$$
which will be repeatedly used in the sequel. 

\subsubsection{Rough scale invariance}
We begin with a proposition which is a rough version of scale invariance. In words it says that the typical height of every subtree grafted on $ \mathcal{T}_{n}$ is at most $n ^{-\alpha+\circ(1)}$. Combined with Proposition \ref{prop:unif}, it is the core of the proof of Proposition \ref{prop:compact}. For a stem $u$, recall the notation $ \mathbf{a}(u)$ from Section \ref{sec:not}.
\begin{proposition} \label{prop:tech1} If $a_{i} \leq i^{-\alpha+\circ(1)}$ for some $\alpha \in (0,1]$, then,  almost surely, 
$$ 
\sup_{u \in \mathcal{G}_{n}} \mathsf{H} \big(  \mathbf{a}(u) \big) \leq n^{-\alpha + \circ(1)}.$$
\end{proposition}
 
\proof We first prove that the longest length of a stem of $ \mathcal{T}_{n}$ is at most $n^{-\alpha+\circ(1)}$. To see this, suppose by contradiction that a stem of length at least $ n^{-\alpha + \varepsilon}$ is present in $ \mathcal{T}_{n}$ for some $ \varepsilon>0$. Provided that $n$ is large enough, since  $a_{i} \leq i^{-\alpha+\circ(1)}$, this stem must be part of a branch $\mathrm b_i$ (of length $a_{i}$) grafted at some time $i \leq n/2$. It thus means that we can find a part of length $n^{-\alpha+ \varepsilon}/2$ of the branch $\mathrm b_{i}$ whose endpoints are exactly at distance $k n^{-\alpha+ \varepsilon}/2$ and $(k+1) n^{-\alpha+ \varepsilon}/2$ for some $k \geq 0$ from the extremity of $\mathrm b_i$ closest to root of $ \mathcal{T}_{i}$ which has not been hit by the grafting process between times $\lfloor n/2 \rfloor +1$ and $n$. For each $k$, such an event has probability at most  
\begin{eqnarray*}\left(1 - \frac{n^{- \alpha + \varepsilon}}{2A_{n}}\right)^{n/2} &{\leq}& \exp(- n^{ \varepsilon + \circ(1)}),  \end{eqnarray*}
since $A_n \leq n^{1-\alpha+\circ(1)}$ because $a_{i} \leq i^{-\alpha + \circ(1)}$ and $\alpha \in (0,1]$.
Summing over all possibilities to choose such a part on some $\mathrm b_i$ for some $i \leq n$, we find that asymptotically the probability that there is a stem of length at least $n^{-\alpha + \varepsilon}$ in $ \mathcal{T}_{n}$ is bounded above by 
 \begin{eqnarray*} \sum_{i \leq n/2} \left(\frac{2a_{i}}{n^{- \alpha + \varepsilon}}+1\right) \exp(-n^{ \varepsilon+\circ(1)})  &=& \exp(-n^{ \varepsilon + \circ(1)}).  \end{eqnarray*}
We easily conclude by an application of Borel--Cantelli  that \begin{eqnarray} \sup_{u \in \mathcal{G}_{n}}   \mathsf{L}( {u}) \leq n^{-\alpha + \circ(1)}.\label{eq:subbrancsmall}  \end{eqnarray}
To deduce from this the proposition, we  
need the following lemma.
\begin{lemma}  \label{lem:expo2} Pick a stem ${u}$ of $ \mathcal{T}_{n}$, then, conditionally on $\mathcal T_n,$  for any $\lambda \geq 0$ such that $\lambda \mathsf{L}({u}) <1$ and $ \lambda a_{i} <1$ for all $i \geq n$, we have
$$ \mathbb{E}\Big[ e^{\lambda  \mathsf{H}( \mathbf{a}(u))} \mid \mathcal{F}_{n} \Big] \leq \exp\left( 2\lambda \bigg( \mathsf{L}({u}) + \sum_{i=n+1}^\infty	\frac{a_{i}^2}{A_{i}}\bigg)\right).$$
\end{lemma}

\proof For $i \geq 1$, let $A_{i}(u) = a_{1}(u) + ... + a_{i}(u)$ and for $p \geq 1$, let $$ \Sigma_{p} = \sum_{i=1}^p \frac{a_{i}(u)^2}{A_{i}(u)}, \qquad \mbox{so that } \Sigma_{\infty} = \mathsf{H}( \mathbf{a}(u)),$$ with the convention that $\frac{a_{i}(u)^2}{A_{i}(u)} =0$ if $a_{i}(u) =0$. Next, let  $\lambda \geq 0$ satisfy the assumptions of the statement.  
For $p \geq n$, since the branch $a_{p+1}$ is grafted on $ \mathcal{T}_{p}(u)$ with probability $A_{p}(u)/A_{p}$, we have
 \begin{eqnarray*} \mathbb{E}\left[e^{\lambda \Sigma_{p+1}} \mid \mathcal{F}_{p}\right]&=& e^{\lambda \Sigma_{p}} \left( \frac{A_{p}-A_{p}(u)}{A_{p}} + \frac{A_{p}(u)}{A_{p}}e^{\lambda \frac{a_{p+1}^2}{A_{p}(u)+a_{p+1}}}\right) \\ & =& e^{\lambda \Sigma_{p}} \left(1+ \frac{A_{p}(u)}{A_{p}}\Big(e^{\lambda \frac{a_{p+1}^2}{A_{p}(u)+a_{p+1}}}-1\Big)\right) \\
 & {\leq} & e^{\lambda \Sigma_{p}} \left( 1 + 2 \lambda \frac{a_{p+1}^2}{A_{p}(u)+a_{p+1}} \times \frac{A_{p}(u)}{A_{p}} \right).
 \end{eqnarray*} 
To go from the second to the third line, we have used that $ \lambda \frac{a_{p+1}^2}{A_{p}(u)+a_{p+1}} \leq \lambda a_{p+1} \leq 1$  and that $e^x-1 \leq 2 x$ for $x \in [0,1]$. 
Besides, since for a fixed $c>0$, the function $x \mapsto x/(x+c)$ is increasing on $(0,\infty)$ and $A_{p}(u) \leq A_{p}$ we have that $\frac{A_{p}(u)}{(A_{p}(u)+a_{p+1})A_p} \leq \frac{1}{A_{p+1}}$, which finally leads to
$$
 \mathbb{E}\left[e^{\lambda \Sigma_{p+1}} \mid \mathcal{F}_{p}\right] \ \leq \ e^{\lambda \Sigma_{p}} \left( 1 + 2 \lambda \frac{a_{p+1}^2}{A_{p+1}}  \right).
$$
Note that we also have  $ \mathbb{E}[e^{\lambda \Sigma_{n}}] = e^{\lambda \mathsf{L}({u})} \leq 1 + 2 \lambda  \mathsf{L}({u})$. So, conditioning in cascades over all integers $p \geq n$, we obtain 
 \begin{eqnarray*}  \mathbb{E}[ e^{\lambda \mathsf{H}( \mathbf{a}(u))}\mid \mathcal{F}_{n}]=  \mathbb{E}[e^{\lambda \Sigma_{\infty}}\mid  \mathcal{F}_{n}] &\leq& (1+ 2 \lambda \mathsf{L}({u}))\prod_{i=n+1}^\infty\left( 1 + 2\lambda \frac{a_{i}^2}{A_{i}} \right) \\ & \leq & \exp\left( 2\lambda \bigg( \mathsf{L}({u}) + \sum_{i=n+1}^\infty	\frac{a_{i}^2}{A_{i}}\bigg)\right).	  \end{eqnarray*}
\endproof
Coming back to the proof of Proposition \ref{prop:tech1}, fix $ \varepsilon>0$ and consider $n_{\eps}$ such that $a_{n}\leq n^{-\alpha+\eps}$ and  $\sum_{n}^\infty \frac{a_{i}^2}{A_{i}} \leq n^{-\alpha + \varepsilon} $ for all $n \geq n_{\eps}$ ($n_{\eps}$ exists by  Lemma \ref{lem:manip} (ii) and since $a_{i} \leq i^{-\alpha+\circ(1)}$). Then, for $m \geq n_{\eps}$, let
$\mathcal{E}_{m}$ denote the event 
$$\sup_{u \in \mathcal{G}_{n}} \mathsf{L}({u})  \leq n^{-\alpha + \varepsilon} \quad \text{ for all } n \geq m.$$ By the first part of the proof, $\mathbb P(\mathcal E_m)$ converges to $1$ as $m \to \infty$. Next, for a fixed $m \geq n_{\eps}$ and all $n\geq m$, using a standard Markov exponential inequality and Lemma \ref{lem:expo2} with $\lambda = n^{\alpha - \varepsilon}$ on the event $\mathcal E_m$, we get  $$ \mathbb{P}\left( \mathsf{H}( \mathbf{a}(u)) \geq n^{-\alpha + 2 \varepsilon} \mid  \mathcal{E}_{m} \right) \leq \frac{ e^{-\lambda n^{- \alpha + 2 \varepsilon}} \mathbb{E}\left[\mathbb{E}\left[e^{\lambda \mathsf{H}( \mathbf{a}(u))} \mathbbm 1_{\mathcal E_m} \mid \mathcal F_n \right]\right]}{\mathbb P(\mathcal E_m)}\leq \frac{e^{-\lambda n^{- \alpha + 2 \varepsilon}+4 \lambda n^{-\alpha + \varepsilon}}}{\mathbb P(\mathcal E_m)} \leq e^{ -n^{\varepsilon+\circ(1)}}.$$
 Since their are $2n-1$ stems in $ \mathcal{T}_{n}$, the Borel--Cantelli lemma  shows that conditionally on $ \mathcal{E}_{m}$ we have $ \sup_{u \in \mathcal{G}_{n}} \mathsf{H}( \mathbf{a}(u)) \leq n^{-\alpha +\circ(1)}$ almost surely. The conclusion follows, since $ \mathbb{P}( \mathcal{E}_{m}) \to 1$ as $m \to \infty$. \endproof 
 
 \begin{remark} When $a_i \leq i^{-\alpha + \circ(1)}$ for some $\alpha > 1$ the statement of this proposition is no longer true. Indeed, in this case the length of the largest stem of $ \mathcal{T}_{n}$ is roughly of order $n^{-1} \gg n^{-\alpha}$. 
 \end{remark}

\subsubsection{Proof of Proposition \ref{prop:compact}}
 
\noindent \textsc{Compactness.} Recall that $ \mathcal{T}_{n}$ and $ \mathcal{T}$ have been built as closed subsets of $ \ell^{1}( \mathbb{R})$. Since the set of non-empty compact subspaces of $ \ell^{1}( \mathbb{R})$ endowed with the Hausdorff distance (denoted here by $\delta_{\mathrm H}$) is complete, it suffices  to show that 
 \begin{eqnarray} \sum_{i\geq 1} \delta_{\mathrm H}( \mathcal{T}_{2^{i+1}}, \mathcal{T}_{2^{i}}) < \infty \quad \mbox{ almost surely} \label{eq:goalcompact}  \end{eqnarray}
 to get the almost sure compactness of $ \mathcal{T}$.
Note that $\delta_{\mathrm H}( \mathcal{T}_{2^{i+1}}, \mathcal{T}_{2^{i}})$ is less than, or equal to, the maximal  height of subtrees $ \mathcal{T}_{2^{i+1}}(u)$ when $u$ runs over $\mathcal{G}_{2^{i}}$ (the subtrees $\mathcal T_n(u), \mathcal T(u)$ are defined in Section \ref{sec:not}). To approximate the heights of these subtrees, we will throw $2^i$  independent uniform points in each of them and take the maximal height attained. Fix $ \varepsilon>0$ and let $n_{\eps}$ be such that $a_n \leq n^{\eps-\alpha}$ for $n \geq n_{\eps}$. For each $m \geq n_{\eps}$, consider the event $ \mathcal{E}'_{m}$ on which \begin{eqnarray*} \label{eq:condi} \sup_{u \in \mathcal{G}_{n}} \mathsf{H}( \mathbf{a}(u)) \leq n^{\eps-\alpha} \quad \text{ for all }n\geq m.  \end{eqnarray*}
By Proposition \ref{prop:tech1}, $ \mathbb{P}( \mathcal{E}'_{m}) \to 1$ as $ m \to \infty$. It thus suffices to work conditionally on $ \mathcal{E}_{m}'$. 

So, fix $i\geq 1$ such that $2^i \geq m$, pick $u \in \mathcal{G}_{2^i}$ and let $H(u)$ denote the height of a random uniform point in $ \mathcal{T}_{2^{i+1}}(u)$. By Proposition \ref{prop:unif} with $\lambda = 2^{i(\alpha - \varepsilon)}$ we have 
 \begin{eqnarray} \label{eq:1} \mathbb{P}\big( H(u) \geq 2^{i(2 \varepsilon - \alpha)}   \mid {\mathcal{E}'_{m}, \mathcal{F}_{2^i}}\big) &\underset{ \mathrm{Markov}}{\leq}&  \frac{\mathbb{E}\left[ e^{2^{i(\alpha - \varepsilon)} H(u)} \mathbbm 1_{ \mathcal{E}_{m}'} \mid \mathcal{F}_{2^i} \right]}{\exp( 2^{i(2 \varepsilon - \alpha)}{ 2^{i( \alpha - \varepsilon)}}) \mathbb P(\mathcal E'_m)} \nonumber \\
  &{\leq }&  \frac{\mathbb{E}\left[ \mathbb{E}\left[  e^{2^{i(\alpha - \varepsilon)} H(u)} \mathbbm 1_{\{\mathsf{H}( \mathbf{a}(u)) \leq 2^{i(\eps-\alpha)} \} } \mid \mathbf{a}(u), \mathcal{F}_{2^i} \right] \mid \mathcal{F}_{2^i} \right]}{\exp(2^{i \varepsilon})\mathbb P(\mathcal E'_m)}   \nonumber \\
 & \underset{ \mathrm{Prop.} \ref{prop:unif}}{\leq}&  \frac{\mathbb{E}\left[ e^{2^{i(\alpha - \varepsilon)} \mathsf{H}( \mathbf{a}(u))}  \mathbbm 1_{\{\mathsf{H}( \mathbf{a}(u)) \leq 2^{i(\eps-\alpha)} \} }\mid \mathcal{F}_{2^i}\right]}{\exp( 2^{i \varepsilon}) \mathbb P(\mathcal E'_m)} \leq \frac{e^1}{\exp( 2^{i \varepsilon})\mathbb P(\mathcal E'_m)}. \label{eq:height2i} \end{eqnarray}
To apply Proposition \ref{prop:unif} in the third line we had to notice that conditionally  on the sequence $ \mathbf{a}(u)$, the tree $ \mathcal{T}(u)$ is constructed from $ \mathbf{a}(u)$ as $ \mathcal{T}$ is constructed from $ \mathbf{a}$. In particular, according to the discussion preceding Proposition \ref{prop:unif}, the height of a uniform point in $ \mathcal{T}_{2^{i+1}}(u)$ is stochastically at most the height of a uniform point in $ \mathcal{T}(u)$, conditionally on  $ \mathbf{a}(u)$.

 We now throw $2^i$ independent uniform points in each of the $2^{i+1} -1$ subtrees $ \mathcal{T}_{2^{i+1}}(u)$, for each $u \in \mathcal{G}_{2^i}$. Let $\mathcal B_i$ denote the event ``the maximal height attained by one of these $(2^{i+1}-1) \cdot 2^i$ uniform points  is at least $2^{i(2 \varepsilon - \alpha)}$". By \eqref{eq:height2i}, conditionally on $ \mathcal{E}_{m}'$, the probability of $\mathcal B_i$  is bounded from above by 
 $$ (2^{i+1}-1)\cdot 2^i \frac{e^1}{ \exp( { 2^{i\varepsilon}})\mathbb P(\mathcal E_m')}.$$
 The last quantity is summable in $i \geq 0$, hence by Borel--Cantelli we conclude that $ \mathcal{B}_{i}$ happens finitely many often, conditionally on $ \mathcal{E}_{m}'$. 

On the other hand, for each $u \in \mathcal{G}_{2^i}$,  the total length of $ \mathcal{T}_{2^{i+1}}(u)$ is at most $A_{2^{i+1}} \leq 2^{i(1-\alpha+\circ(1))}$. Hence when we throw independently $2^{i}$ uniform points in this subtree, the probability that none of these points is at distance less than $2^{i(2 \varepsilon - \alpha)}$ of the maximal height is at most 
 $$  \left( 1- \frac{2^{i(2 \varepsilon - \alpha)}}{A_{2^{i+1}}}\right)^{2^i}\leq  \exp\left(- 2 ^i\frac{2^{i(2 \varepsilon - \alpha)}}{2^{i(1-\alpha+\circ(1))}}\right) = \exp(-2^{i(2 \varepsilon +\circ(1))}).$$ 
Even after multiplying the right-hand side by $2^{i+1}-1$ the series is still summable, and so after another application of the Borel--Cantelli lemma, we can gather the last two results to deduce that almost surely (conditionally on $ \mathcal{E}_{m}'$) for $i$ large enough the heights of all subtrees $ \mathcal{T}_{2^{i+1}}(u)$, $u \in \mathcal{G}_{2 ^i}$ is at most $2 \cdot 2^{i ( 2 \varepsilon - \alpha)}$.  Letting $m \rightarrow \infty$, this readily leads to \eqref{eq:goalcompact}. 

\medskip

\noindent \textsc{Upper bound on the Hausdorff dimension.} All the assertions in this paragraph hold almost surely. From the previous discussion, we deduce that conditionally on $ \mathcal{E}_{m}'$ the diameter of the trees $ \mathcal{T}(u)$ for $u \in \mathcal{G}_{2^i}$ is at most $2^{i(3 \varepsilon -\alpha)}$ for all $i$ large enough. For those integers $i$, we thus obtain a covering of $ \mathcal{T}$  made of $2^{i+1}-1$ balls of diameter $2^{i(4 \varepsilon -\alpha)}$. This immediately implies that $ \mathrm{dim_{H}}( \mathcal{T}) \leq 1/( \alpha - 4 \varepsilon)$. Since $ \varepsilon>0$ was arbitrary and $ \mathbb{P}( \mathcal{E}_{m}') \to 1$, we indeed proved that $ \mathrm{dim_{H}}( \mathcal{T}) \leq 1/\alpha$ a.s. 
\endproof

\subsection{Lower bound via $\mu$}

Together with Proposition \ref{prop:compact} and the fact $  \mathrm{dim_{H}}(\mathcal{T}) \geq 1$, the following result implies Theorem \ref{thm:main}.

\begin{proposition} \label{prop:lowerbound}
Assume that $\mathsf{H}( \mathbf{a})<\infty$ and  $A_n \geq n^{1-\alpha+\circ(1)}$ for $ \alpha \in (0,1)$. Then, the Hausdorff dimension of $ \mathcal{T}$ is at least $\alpha^{-1}$ almost surely.
\end{proposition}

Note that this result also applies to cases where we do not know if the tree $\mathcal T$ is compact. E.g. the two hypotheses hold when $a_i=\ln(i)^{-\gamma}$ for some $\gamma>1$, for all $\alpha \in (0,1]$. In this case the Hausdorff dimension of the tree is therefore infinite a.s.  

\begin{remark} \label{rem:mieux} When $\mathsf{H}( \mathbf{a})<\infty$ and $A_n \rightarrow \infty$,  our proof below can easily be adapted to show that the Hausdorff dimension of $\mathsf{Leaves}(\mathcal{T})$ is at least 1 almost surely.
\end{remark}

The rest of this section is devoted to the proof of Proposition \ref{prop:lowerbound}. Our approach relies on Frostman's theory and the existence of the measure $\mu$, the weak limit of the uniform measures $\mu_n$ which exists when $ \mathsf{H}( \mathbf{a}) < \infty$ by Theorem \ref{thm:measures}. More precisely, we know by a result of Frostman \cite[Theorem 4.13]{Falc03}, that 
$$
\int_{\mathcal T \times \mathcal T} \frac{\mu(\mathrm d x)\mu(\mathrm d y)}{\left(\delta(x,y)\right)^{\gamma}}<+\infty \ \Rightarrow \ \mathrm{dim}_{\mathrm H}(\mathcal T) \geq \gamma
$$ 
(we recall that $\delta$ denotes the distance on $\mathcal T$).
Hence, given $\mathcal T$, consider two points picked uniformly and independently at random according to the measure $\mu$, and let  $D$ denote their distance in $\mathcal T$. Clearly,
$$
\mathbb E\left[D^{-\gamma} \right]=\mathbb E \left[ \int_{\mathcal T \times \mathcal T} \frac{\mu(\mathrm d x)\mu(\mathrm d y)}{(\delta(x,y))^{\gamma}}\right],
$$
from which we deduce that it is sufficient to prove that $\mathbb  E\big[D^{-\gamma} \big]<\infty$ for all $\gamma \in (0,\alpha^{-1})$ to get the desired lower bound. This will be implied by the following lemma: 

\begin{lemma}
\label{lem:majoDn}
Under the conditions of Proposition \ref{prop:lowerbound}, for all $\eps>0$, $\exists c_{\alpha,\eps}>0$ such that for all $r \in (0,1]$,
\begin{equation*}
\label{eq:majoDn}
\mathbb P\left(D \leq r \right) \leq c_{\alpha,\eps} r^{\frac{1}{\alpha}-\eps}.
\end{equation*}
Consequently, $\mathbb E\big[D^{-\gamma} \big]<\infty$ for all $\gamma \in (0,\alpha^{-1})$.
\end{lemma}

To prove the last lemma we will compute exactly the (annealed) law of $D$ in a similar fashion  we computed the exact law of the height of a random point sampled according to $\mu$. We then proceed to the proof of Lemma \ref{lem:majoDn}. 

\subsubsection{Description of the law of the two-point function}
\label{sec:Dinfty}

\begin{lemma} 
\label{lem:Dinfini}
Let $U_{i},V_{i},V'_{i}, i \geq 1$ be random variables independent and uniform on $[0,1]$. The distribution of $D$ is given by
$$
\mathbb E \left[f(D)\right]=\sum_{k=1}^{\infty} \Bigg[\left( \frac{a_{k}}{A_{k}} \right)^2 \prod_{j=k+1}^{\infty} \bigg(1 - \left(\frac{a_{j}}{A_{j}}\right)^2\bigg)\Bigg] \mathbb E \left[f\left(a_{k}|V_{k}-V_{k}'|+ \sum_{i=k+1}^\infty a_{i}V_{i} \mathbbm{1}_{\left\{U_{i} \leq 2\frac{a_{i}}{A_{i}+a_i} \right\}} \right) \right]
$$
for all measurable positive functions $f$.
\end{lemma}

\proof 
Let $ n \geq 2$ and conditionally on $ \mathcal{T}$ consider two points $Y_{n}^{(1)}$ and  $Y_{n}^{(2)} \in \mathcal{T}_{n}$  independent and distributed according to $\mu_{n}$. We let $D_n$ denote their distance.
\begin{itemize}
\item With probability $(1 - \frac{a_{n}}{A_{n}})^2$ these two points belong to $ \mathcal{T}_{n-1}$ and conditionally on this event they are independent, uniform on $ \mathcal{T}_{n-1}$. On this event we thus have $D_{n} \overset{\mathrm{(d)}}=D_{n-1}$.
\item With probability $2 (1 - \frac{a_{n}}{A_{n}}) ( \frac{a_{n}}{A_{n}} )$ only one of these points belongs to the $n$th branch. Conditionally on this event, the point in question is uniformly distributed on the last branch and the remaining point is independent and uniform on $ \mathcal{T}_{n-1}$. Moreover the projection of these two points onto $ \mathcal{T}_{n-1}$ yields a pair of independent points uniformly distributed over $ \mathcal{T}_{n-1}$. On this event we thus have $D_{n} \overset{\mathrm{(d)}}= D_{n-1} + a_{n} V_{n}$ where in the right side, $V_{n}$ is uniform on $(0,1)$ and independent of $D_{n-1}$. 
\item Finally, with probability $( \frac{a_{n}}{A_{n}} )^2$ these two points belong to the $n$th branch. Conditionally on this event they are uniform, independent on this branch,
and thus we can write $D_{n} = a_{n}|V_{n}-V_{n}'|$ where $V_{n}$ and $V_{n}'$ are independent and both uniform on $(0,1)$.
\end{itemize}
Noticing that for $n \geq 2$  $$ \frac{2 (1 - \frac{a_{n}}{A_{n}}) ( \frac{a_{n}}{A_{n}} )}{ 1 - ( \frac{a_{n}}{A_{n}} )^2} =  \frac{2 a_{n}}{A_{n}+a_n},$$  it follows from the previous discussion that the law of $D_{n}$ is described as follows:
 \begin{eqnarray*} \label{eq:loiDn}\begin{array}{c} \mbox{for\ } k \in \{1,2, ... , n\}\ \  \mbox{with probability }  \displaystyle \left( \frac{a_{k}}{A_{k}} \right)^2 \prod_{i=k+1}^n\bigg(1 - \left(\frac{a_{i}}{A_{i}}\right)^2\bigg) \\ 
 \mbox{ we have }  \displaystyle D_{n} = a_{k}|V_{k}-V_{k}'|+ \sum_{i=k+1}^n a_{i}V_{i} \mathbbm{1}_{\left\{U_{i} \leq 2\frac{a_{i}}{A_{i}+a_i} \right\}}, \end{array} \end{eqnarray*} where the variables $U_{i},V_{i},V'_{i}, 1 \leq i \leq n$ are all independent and uniform on $[0,1]$ (we use the convention that the sum over the empty set is $0$, whereas the product over the empty set is $1$). From Theorem \ref{thm:measures}, we get that
$
D_n \to D$ in distribution so that  passing to the limit, we get a similar description of the law of $D$. In this last step, it is crucial that the series $\sum_k (\frac{a_k}{A_k})^2$ converges to ensure that $\mathbb P(D=\infty)=0$. We check in Lemma \ref{lem:manip} that such a series is always convergent. \endproof

\subsubsection{Proof of Lemma \ref{lem:majoDn}}
Fix $\eps \in (0,1)$ and let $r \in (0,1]$. By Lemma \ref{lem:Dinfini} we have 
  \begin{eqnarray*} 
  && \mathbb{P}(D \leq r )  =  \sum_{k=1}^\infty \Bigg[\left( \frac{a_{k}}{A_{k}} \right)^2\prod_{j=k+1}^\infty \bigg(1- \left( \frac{a_{j}}{A_{j}}\right)^2\bigg)\Bigg] \mathbb{P}\bigg(  a_{k}|V_{k}-V_{k}'|+ \sum_{i=k+1}^\infty a_{i}V_{i} \mathbbm{1}_{E_i} \leq r \bigg) \\ 
  & \leq &  \sum_{k= \lfloor r^{-\frac{1}{\alpha}+\frac{\eps}{2}}\rfloor +1}^{+\infty} \left(\frac{a_k}{A_k} \right)^2 \mathbb{P}\left( a_{k}|V_{k}-V_{k}'| \leq r \right) \\
  &+& \sum_{k=1}^{ \lfloor r^{-\frac{1}{\alpha}+\frac{\eps}{2}}\rfloor} \left( \frac{a_{k}}{A_{k}} \right)^2  \mathbb{P}\left( a_{k}|V_{k}-V_{k}'| \leq r \right)\prod_{i=k+1}^{ \lfloor r^{-\frac{1}{\alpha}+\frac{\eps}{2}}\rfloor} \mathbb{P}\left( a_{i}V_{i} \mathbbm{1}_{E_i} \leq r \right)  \bigg(1- \left( \frac{a_{i}}{A_{i}}\right)^2\bigg),
  \end{eqnarray*}
where we have set $E_i=\Big\{U_{i} \leq 2\frac{a_{i}}{A_{i}+a_i} \Big\}$ to improve the presentation.
Then, note that
\begin{eqnarray*}
\mathbb{P}\left(a_{i}V_{i} \mathbbm{1}_{E_i} \leq r \right) \ \leq \ 1-\frac{2a_i}{A_{i}+a_i}  +  \frac{2a_i}{A_{i}+a_i} \times \frac{r}{a_i} \ \leq \  \frac{A^2_{i-1}}{A^2_i} \times \frac{A_{i}^2}{A_i^2-a_i^2} \times \left(1+\frac{2r}{A_{i-1}} \right),
\end{eqnarray*} 
which leads us to
\begin{eqnarray*} 
\prod_{i=k+1}^{\lfloor r^{-\frac{1}{\alpha}+\frac{\eps}{2}}\rfloor}  \mathbb{P}\left(a_{i}V_{i} \mathbbm{1}_{E_i} \leq r \right)  \bigg(1- \left( \frac{a_{i}}{A_{i}}\right)^2\bigg)
\leq 
\prod_{i=k+1}^{\lfloor r^{-\frac{1}{\alpha}+\frac{\eps}{2}}\rfloor} \frac{A^2_{i-1}}{A^2_i} \times   \prod_{i=k+1}^{\lfloor r^{-\frac{1}{\alpha}+\frac{\eps}{2}}\rfloor} \left(1+\frac{2r}{A_{i-1}} \right).
\end{eqnarray*}
But the second product in the right-hand side is bounded from above by a constant independent of $k$ and $r \in (0,1]$. Indeed, using that $\ln(1+x)\leq x$ for positive $x$, we get that
$$
\prod_{i=k+1}^{\lfloor r^{-\frac{1}{\alpha}+\frac{\eps}{2}}\rfloor} \left(1+\frac{2r}{A_{i-1}} \right) \leq \exp\Bigg( 2r\sum_{i=k+1}^{\lfloor r^{-\frac{1}{\alpha}+\frac{\eps}{2}}\rfloor} \frac{1}{A_{i-1}} \Bigg) \leq \exp\left( 2r (r^{-\frac{1}{\alpha}+\frac{\eps}{2}})^{\alpha + \circ(1)}\right),
$$ 
where we have used the assumption on the lower bound of $A_n$ for the second inequality (here the notation $\circ$ refers to the convergence of $r$ towards $0$).
Finally, we have proved the existence of a finite constant $C$ independent of $r \in (0,1]$ such that 
$$
\mathbb P(D \leq r) \leq \sum_{k= \lfloor r^{-\frac{1}{\alpha}+\frac{\eps}{2}}\rfloor+1}^{+\infty} \left(\frac{a_k}{A_k} \right)^2 \times \frac{2r}{a_k}+ C\sum_{k=1}^{ \lfloor r^{-\frac{1}{\alpha}+\frac{\eps}{2}}\rfloor} \left(\frac{a_k}{A_k} \right)^2 \times \frac{2r}{a_k}\times  \frac{A_k^2}{A^2_{\lfloor r^{-\frac{1}{\alpha}+\frac{\eps}{2}}\rfloor}}. 
$$
By Lemma \ref{lem:manip} (iii), the first sum in the right-hand side is at most $ r^{\frac{1}{\alpha}-\frac{(1-\alpha)\eps}{2} +\circ(1)}$. So we finally get,
\begin{eqnarray*}
\mathbb P(D \leq r)&\leq &r^{\frac{1}{\alpha}-\frac{(1-\alpha)\eps}{2} +\circ(1)} +  \frac{2rC}{A_{\lfloor r^{-\frac{1}{\alpha}+\frac{\eps}{2}}\rfloor}}\\
&\leq & r^{\frac{1}{\alpha}-\frac{(1-\alpha)\eps}{2} +\circ(1)}.
\end{eqnarray*}

\section{Finite length case}
The goal of this section is to prove Theorem \ref{thm:maina>1}. As in the previous section, we will first prove the compactness and the upper bound of the Hausdorff dimension, which hold in a more general (and even deterministic) setting than that of Theorem \ref{thm:maina>1}. The lower bound on the dimension is more technical than in the previous section and requires the construction of a new measure supported by the leaves of $ \mathcal{T}$.

\subsection{Deterministic results in the finite length case}
The following proposition does not depend on the fact that the new branches are grafted uniformly on the pre-existing tree, but just on the asymptotic behavior of the sequence $(a_i,i\geq 1)$. So, in this subsection, and only in this subsection,  $\mathcal{T}$ designs the completion of a tree built by grafting the branches $\mathrm b_i$ of lengths $a_{i}$ iteratively, without any explicit rules on where the branches are glued. We denote by $ \mathsf{Leaves}( \mathcal{T})$ the set of leaves of $ \mathcal{T}$. 

\begin{proposition}[] \label{prop:finitelength} If $\sum_{i=1}^\infty a_{i} < \infty$, the tree $ \mathcal{T}$ is compact and of Hausdorff dimension $1$. Moreover,
$$
\mathrm{dim}_{\mathrm H}(\mathsf{Leaves} (\mathcal T)) \leq \gamma \quad \mbox{ as soon as }  \quad \sum_{i=1}^\infty a_{i}^\gamma < \infty .
$$
\end{proposition} 
\proof We start with the proof of the upper bound of the Hausdorff dimension of the leaves and assume that $\sum_{i\geq 1}a^{\gamma}_i<\infty$ for some $\gamma \leq 1$. Since the set of leaves of $\mathcal T^*$ is at most countable, its Hausdorff dimension is 0. To get the expected upper bound, we thus only need to get an upper bound for the Hausdorff dimension of $\mathcal T \backslash \mathcal T^*$. 

In that aim, fix $\eps>0$ and let $n_{\eps}$ be such that $\sum_{i > n_{\eps}} a_i \leq \eps.$ Consider then the decomposition of $\mathcal T \backslash \mathcal T_{n_{\eps}}$ 
into connected components and note that the set of closures of these components forms a (at most) countable set of closed subtrees of $\mathcal T$, that covers $\mathcal T\backslash \mathcal T^*$. The intersection of each of these subtrees with  $\mathcal T_{n_{\eps}}$ is reduced to a unique point, the root of the subtree (different subtrees may have the same root -- recall that we have no explicit rule of gluing). We denote by $\mathcal R_{\eps}$ this set of roots, and, for all $\mathrm r \in \mathcal R_{\eps}$, by $\mathcal T_{n_{\eps}}^{(\mathrm r)}$ the union of subtrees descending from it, which is also a tree.  We then let $\mathcal I_{\mathrm r}$ be the set of integers $i$ such that the segment $\mathrm b_i$  belongs to the subtree $\mathcal T_{n_{\eps}}^{(\mathrm r)}$. Clearly, this subtree has a diameter at most $\sum_{i \in \mathcal I_{\mathrm r}} a_i$ which is itself at most $\eps$, by definition of $n_{\eps}$.

The collection of subtrees $\mathcal T_{n_{\eps}}^{(\mathrm r)}$, $\mathrm r \in \mathcal R_{\eps}$ therefore forms an at most countable covering of $\mathcal T \backslash \mathcal T^*$ with sets of diameter less than $\eps$.
We have
$$
\sum_{\mathrm r \in \mathcal R_{\eps}} \bigg(\sum_{i \in \mathcal I_{\mathrm r}} a_i\bigg)^{\gamma} \leq \sum_{\mathrm r \in \mathcal R_{\eps}} \sum_{i \in \mathcal I_{\mathrm r}} a_i^{\gamma} \leq \sum_{i\geq 1} a_i^{\gamma}<\infty,
$$
where the first inequality holds since $\gamma \leq 1$ and the second since the sets $\mathcal I_{\mathrm r}, \mathrm r \in \mathcal R_{\eps}$ are disjoint. Hence the $\gamma-$dimensional Hausdorff measure of $\mathcal T \backslash \mathcal T^*$ is finite and its Hausdorff dimension is at most $\gamma$ (almost surely).

We now turn to the compactness  of $\mathcal T$ under the sole assumption $\sum_{i \geq 1} a_{i}< \infty$. We consider $\eps>0$ and use the notation introduced above. The tree $\mathcal T_{n_{\eps}}$ is clearly compact and we let $B(x_n,\eps)$, $n \leq  N_{\eps}$ be a finite collection of open balls of radius $\eps$ that covers it. Besides, as noticed above, all $x \in \mathcal T \backslash \mathcal T_{n_{\eps}}$ is at distance at most $\eps$ from an element of $\mathcal R_{\eps}$. Consequently the collection of open balls $B(x_n,2\eps)$, $n \leq  N_{\eps}$ of radius $2\eps$ covers $\mathcal T$. Hence $ \mathcal{T}$ is pre-compact and thus compact by completeness. \endproof

\subsection{Lower bound for the Hausdorff dimension of the leaves}

In this section we assume the existence of $\alpha>1$ such that
$$ (\mathrm D_{\alpha}) \qquad a_{i} \leq i^{-\alpha+\circ(1)} \quad \mbox{and} \quad a_{i}+ a_{i+1}+ ...+ a_{2i} = i^{1- \alpha + \circ(1)}.$$
In particular,  by Proposition \ref{prop:finitelength}, the tree $ \mathcal{T}$ is compact and the Hausdorff dimension of its set of leaves is bounded above by $1/\alpha$ (almost surely). The following result is the complement to obtain the statement of Theorem \ref{thm:maina>1}.
\begin{proposition} \label{prop:lowbound}Under $(\mathrm D_{\alpha})$, almost surely,
$$\mathrm{dim}_{\mathrm H}(\mathsf{Leaves} (\mathcal T)) \geq 1/\alpha.$$
\end{proposition}
To get this lower bound,  we will show that for any $ \varepsilon>0$ we can construct, with a probability at least $1 - \varepsilon$, a (random) probability measure $\pi$ supported by the set of leaves of $ \mathcal{T}$ such that  for every $x \in \mathcal{T}$
\begin{eqnarray} 
\label{eq:lowerboundgoal}
\limsup_{r \to 0} \frac{\pi \big(B(x,r)\big)}{r^{\frac{1}{\alpha}-\varepsilon}} &=& 0 ,
\end{eqnarray}
where $B(x,r)$ denotes the open ball in $\mathcal T$ of radius $r$ centered at $x$.
By standard results on Hausdorff dimensions (see e.g. \cite[Proposition 4.9]{Falc03}), this will entail that $ \mathrm{dim}_{\mathrm H}( \mathsf{Leaves}( \mathcal{T})) \geq \alpha^{-1} - \varepsilon$ with probability at least $1- \varepsilon$. (Proposition 4.9 in \cite{Falc03} is stated for subsets of $\mathbb R^n$, but, clearly, its proof also holds for any metric space.) Since $ \varepsilon>0$ is arbitrary, this will prove Proposition \ref{prop:lowbound}. \medskip

From now on, $ \varepsilon \in (0,1/\alpha)$ is fixed.  Rather than tempting to construct a ``uniform" measure on the leaves of $ \mathcal{T}$, the support of $\pi$ will be a strict subset of $ \mathsf{Leaves}( \mathcal{T})$. To construct this measure, we need some more notation. 

\medskip

\noindent \textbf{Subsets of good branches}. For $i \geq 1$, we say that the branch $ \mathsf{b}_{i}$, of length $a_i$, is ``good'' if  $i^{-\alpha- \varepsilon} \leq a_{i} $. In other words, a good branch is not too small when it appears (it cannot be greater than $i^{-\alpha + \varepsilon}$ eventually according to $(\mathrm D_{\alpha}))$. For $n \geq 1$,  let $$G_{n} = \{ i \in [[n, 2n]] : \mathsf{b}_{i} \mbox{ is good}\} \quad \text{and} \quad \ell_{n} = \sum_{i \in G_{n}} a_{i},$$ $\ell_n$ being the total length of good branches of index between $n$ and $2n$. It is easy to see  that under assumption $(\mathrm D_{\alpha})$
 \begin{eqnarray} \label{eq:goodbranches} \# G_{n} = n^{1+\circ(1)}  \quad \mbox{and}\quad \ell_{n} = n^{1-\alpha+\circ(1)}.  \end{eqnarray}
Let now $1= n_1<n_2<n_3\ldots$ be integers such that $n_{k+1} > 2n_{k}$ for all $k \geq 1$. Later we will need to do some additional assumptions on the integers $n_k$'s ensuring that they grow sufficiently fast, but for the moment we stay on this. For $\mathsf{b}_{i} , \mathsf{b}_{j}$ two good branches with indices $1\leq j<i$, we write $ \mathsf{b}_{i} \to \mathsf{b}_{j}$ if $\mathsf b_i$ is directly grafted on $\mathsf b_j$. We let $  \mathcal{B}_{1}  = \mathsf{b}_{1}$ and for $k\geq 2$ we define recursively  the subsets $ \mathcal{B}_{k}$ of $ \mathcal{T}$, by deciding that $\mathcal B_k$ is made of the good branches $ \mathsf{b}_{i_{k}}$, $n_k\leq i_k\leq 2n_k$ that are grafted on  (good) branches of $\mathcal B_{k-1}$. This leads to branches of the form 
 \begin{eqnarray*} 
 \mathsf{b}_{i_{k}} \to \mathsf{b}_{i_{k-1}} \to ...\to \mathsf{b}_{i_{2}} \to  \mathsf{b}_{1}  \quad \mbox{ with } \quad n_{\ell} \leq i_{\ell} \leq 2 n_{\ell} \ \mbox{ for every }2 \leq \ell \leq k.  
 \end{eqnarray*}
Note that the sets $ \mathcal{B}_{k},k\geq 1$ may be empty. Slightly changing the notation introduced in Section \ref{sec:not}, we let
$$
\mathcal T(\mathsf b_i)=\big\{x\in \mathcal T: [x]_i \in \mathsf b_i\big\}
$$
denote the subtree descending from $\mathsf b_i$ and
$$
\mathcal T(\mathcal B_k)=\bigcup_{i:\mathsf b_i \in \mathcal B_k} \mathcal T(\mathsf b_i).
$$
Remark that $\mathcal T(\mathcal B_{k+1}) \subset \mathcal T(\mathcal B_{k})$ for all $k \geq 1$. Conditionally on the event $\{ \mathcal{B}_{k} \ne \varnothing, \forall k \geq 1\}$, let now $\pi_{k}$ denote the normalized length measure on $ \mathcal{B}_{k}$. We will see later, choosing the $n_k$'s adequately, that the probability of this  event  can be made arbitrary close to $1$ and that \emph{the measure $\pi$ will  be obtained as a (subsequential) limit of $(\pi_{k})_{k \geq 1}$}. Remark  that conditionally on $\{ \mathcal{B}_{k} \ne \varnothing, \forall k \geq 1\}$, the family $(\pi_{k})_{k\geq1}$ is a sequence of probability measures on a compact space, hence it admits at least one subsequential limit. 
We begin with a simple lemma. 
\begin{lemma} \label{lem:leaves} Almost surely, conditionally on $\{ \mathcal{B}_{k} \ne \varnothing, \forall k \geq 1\}$ (and provided that this event has a positive probability) any subsequential limit $\varpi$ of  $(\pi_k)_{k\geq 0}$ is supported by $\bigcap_{k \geq 1} \mathcal{T}( \mathcal{B}_{k})$, which is included in the set of leaves of $ \mathcal{T}$.
\end{lemma}
 \proof Clearly,  $ \delta (  \mathcal{T}(\mathcal{B}_{k+1}), \mathcal{T}(\mathcal{B}_{k})^c)>0$ almost surely for all $k \geq 1$. Hence we can find an open set $ \mathcal{O}_{k}$ containing $ \mathcal{T}(\mathcal{B}_{k})^c$ such that $\pi_{j}( \mathcal{O}_{k}) =0$ for all $j \geq k+1$ and all $k$, a.s. By the Portmanteau theorem,  it follows that a.s. for any subsequential limit $\varpi$ of  $(\pi_k)_{k\geq 0}$, $\varpi( \mathcal{O}_{k}) =0$ for all $k$ and so
   \begin{eqnarray*} \label{eq:support} \mathrm{Supp}(\varpi) \subset \bigcap_{k \geq 1} \mathcal{T}( \mathcal{B}_{k}). \end{eqnarray*}
   Since $\mathcal T({\mathcal B_k}) \subset \mathcal T \backslash \mathcal T_{n_k-1}$ for all $k$, the right-hand side is a subset of  $ \mathcal{T} \backslash \mathcal T^*$. \endproof 

\medskip

\subsubsection{Lengths estimates.} 
Before embarking into the proof of Proposition \ref{prop:lowbound}, we have to set up some estimates on the total length of descendants in $\mathcal B_{k+1}$ of a given subset of $\mathcal B_k$ and also to check that the distance between most branches of $\mathcal B_k$ is not too small provided that the sequence $(n_k)$ grows sufficiently fast.  This is the goal of this subsection. Once this will be done, we will see in the next subsection how to use this to show that when the sequence $(n_k)$ grows sufficiently fast, the number of branches composing $\mathcal B_k$ is roughly of order $n_k$ whereas their lengths are of order $n_k^{-\alpha}$. This is a first hint that any subsequential limit of  $(\mathcal \pi_k)$ should satisfy (\ref{eq:lowerboundgoal}). Of course, we will need to control our  approximations and the material to do that is developed here. 
We start with some estimates of the total length of good branches indexed by $G_{n}$ that are grafted on a given subset of $\mathcal T_{n-1}$, $n \geq 1$.

\begin{lemma} \label{lemma:tech} Let $n \geq 2$ and consider a subset $S \subset \mathcal{T}_{n-1}$ measurable with respect to $ \mathcal{F}_{n-1}$. Denote by $ \mathcal{X}$ the total length of the branches indexed by $G_{n}$ that are (directly) grafted on $S$. 
\begin{enumerate}
\item[$\mathrm{(i)}$] Then for every $ \eta \in (0,1)$ we have
$$ \mathbb{P}\left( \left| \mathcal{X} - \frac{\ell_{n} |S|}{A_{\infty}} \right| \geq \eta \frac{ \ell_{n}|S|}{A_{\infty}}\right) \leq  \frac{n^{- c +\circ(1)}}{|S|\eta^2}, \quad \mbox{with } c = 1 \wedge(\alpha-1)>0.$$
\item[$\mathrm{(ii)}$] Fix $\delta>0$ and $m\in \mathbb N$. Then, for all $n$ large enough and then for all subsets $S$ such that $ |S| \geq n^{-1+\delta}$,
$$
\mathbb E\left[\mathcal X^m\right] \leq C_m (|S| \ell_n)^m,
$$
where $C_m$  depends only on $m$.
\end{enumerate}
\end{lemma}

\proof  By construction, the random variable $ \mathcal{X}$ can be written as follows: 
$$  \mathcal{X} = \sum_{i \in G_{n}}  a_{i} \mathbbm{1}_{ \left\{U_{i} \leq \frac{|S|}{A_{i-1}} \right\}},$$ where $(U_{i})_{i\geq 1}$ is a sequence of independent random variables uniformly distributed on $(0,1)$.  In particular, $\mathbb E\left[\mathcal X \right]=\sum_{i\in G_n}  \frac{a_i|S|}{A_{i-1}}.$

\noindent (i) Consider temporarily the variable $ \tilde{ \mathcal{X}} = \sum_{i \in G_{n}}  a_{i}\mathbbm{1}_{ \left\{U_{i} \leq \frac{|S|}{A_{\infty}} \right\}}$ instead of $X$. Clearly,  $ \mathbb{E}\big[ \tilde{ \mathcal{X}}\big] = \ell_{n} |S|/A_{\infty}$ and 
 \begin{eqnarray*}  \label{eq:besoin1} 
 \mathrm{Var}\big( \tilde{ \mathcal{X}}\big) \ = \ \sum_{i \in G_{n}} a_i^2 \mathrm{Var}\left(\mathbbm{1}_{ \left\{U_{i} \leq \frac{|S|}{A_{\infty}} \right\} }\right) 
\  = \ \sum_{i \in G_{n}} a_i^2 \left(\frac{|S|}{A_{\infty}}\right)\left(1- \frac{|S|}{A_{\infty}}\right) 
\ \underset{(\mathrm D_{\alpha})}{\leq } \ |S| n^{1-2\alpha+\circ(1)}.  \end{eqnarray*}
On the other hand,  $A_{\infty}-A_{n} = n^{1-\alpha+\circ(1)}$, again by $(\mathrm D_{\alpha})$, and so $$\mathbb{E}\left[\big| \mathcal{X}-\tilde{ \mathcal{X}}\big|\right] \ = \ \sum_{i \in G_{n}} a_{i} \frac{|S|}{A_{\infty}}\frac{(A_{\infty}-A_{i-1})}{A_{i-1}} \ = \ n^{1-\alpha+\circ(1)}\ell_{n}|S|.$$
This leads to
 \begin{eqnarray*} 
 \mathbb{P}\left( \left| \mathcal{X} - \frac{\ell_{n} |S|}{A_{\infty}} \right| \geq 2 \eta \frac{ \ell_{n}|S|}{A_{\infty}}\right) &\leq& \mathbb{P}\left( \left| \tilde{\mathcal{X}} - \frac{\ell_{n} |S|}{A_{\infty}} \right| \geq \eta \frac{ \ell_{n}|S|}{A_{\infty}}\right) + \mathbb{P}\left( \big| \mathcal{X} - \tilde{ \mathcal{X}} \big| \geq \eta \frac{ \ell_{n}|S|}{A_{\infty}}\right)\\ 
 & \leq & \frac{\mathrm{Var}\big( \tilde{\mathcal{X}}\big)}{\eta^2 \ell_{n}^2 |S|^2/A_{\infty}^2} + \frac{ \mathbb{E}\big[ \big|\mathcal{X} -\tilde{ \mathcal{X}}\big|\big]}{ \eta \ell_{n}|S|/A_{\infty}} \\ 
  & \leq & \frac{n^{-1+\circ(1)}}{ |S|\eta^2} + \frac{n^{1-\alpha+\circ(1)}}{\eta}.
 \end{eqnarray*}
 
\noindent (ii) Next, let $i_1,\ldots,i_{\# G_n}$ denote the indices of integers $i \in G_n$. We have for all integers $m \geq 1$,
\begin{eqnarray*}
\mathbb E\left[\mathcal X^{m}\right]&=&  \sum_{\tiny{\begin{array}{c}n_{i_1},\ldots,n_{i_{\# G_n}}: \\n_{i_1}+\ldots+n_{i_{\#G_n}}=m\end{array}}}  \left(\begin{array}{c} m \\ n_{i_1},  \ldots, n_{i_{\# G_n}} \end{array}\right) \prod_{j=1}^{\# G_n} a_{i_j}^{n_{i_j}} \mathbb E\left[\left (\mathbbm 1_{\big\{U_{i_j} \leq \frac{|S|}{A_{i_j-1}}\big\}}\right )^{n_{i_j}}\right] \\
&\leq & m!  \sum_{\tiny{\begin{array}{c}n_{i_1},\ldots,n_{i_{\# G_n}}: \\n_{i_1}+\ldots+n_{i_{\#G_n}}=m\end{array}}} \left(\frac{|S|}{A_1}\right)^{\#\{j:n_{i_j} \geq 1\}}\prod_{j=1}^{\# G_n} a_{i_j}^{n_{i_j}},
\end{eqnarray*}
where we have simply bounded the multinomial term by $m!$. Observe that for every $\# G_n$-tuple involved in the sum, by ($\mathrm D_{\alpha}$), $$\prod_{j=1}^{\# G_n} a_{i_j}^{n_{i_j}}\leq  n^{-m(\alpha  + \circ(1))}.$$ Then, by grouping the $\# G_n$-tuples according to the number of non-zero terms they contain, we get the existence of a constant $c_m$ depending only on $m$ such that
\begin{equation*}
\mathbb E\left[\mathcal X^{m}\right] \leq  m! \sum_{\tiny{\begin{array}{c}n_{i_1},\ldots,n_{i_{\# G_n}} \in \{0,1\}: \\n_{i_1}+\ldots+n_{i_{\#G_n}}=m\end{array}}} \left(\frac{|S|}{A_{1}}\right)^m\prod_{j=1}^{\# G_n} a_{i_j}^{n_{i_j}} + 
c_m \sum_{p=1}^{(m-1)\wedge \# G_n}  {\# G_n\choose p}  |S|^p n^{-m(\alpha  + \circ(1))}.
\end{equation*}
Note that the first term in the right-hand side may be null (if $\# G_n<m$) and is anyway always at most $(A_{1}^{-1}|S| \ell_n)^m$.
Now, noticing that ${\# G_n\choose p} \leq (\# G_n)^{p}$ and using that $|S| \geq n^{-1+\delta}$, we see by (\ref{eq:goodbranches}) that 
$$
 {\# G_n\choose p} |S|^p n^{-m(\alpha  + \circ(1))} \leq (|S| \ell_n)^m,
$$
provided that $n$ is large enough, independently of $p,|S|$. This is sufficient to conclude. \endproof

 \begin{corollary} 
 \label{cor:main}
 There exists a function  $ f : \mathbb{N} \to \mathbb{N}$ with $f(n) > 2n$  for all $n \geq 1$, such that if the sequence $(n_k)_{k \geq 1}$ satisfies $n_{k+1} \geq f(n_{k})$ for all $k \geq 1$, then with probability at least $1- \varepsilon$, 
\begin{equation}
\label{eq:taille} 
  | \mathcal{T}( \mathsf{b}_{i}) \cap \mathcal{B}_{k+1}| \in \left[(1-2^{-k}) \frac{ a_{i} \ell_{n_{k+1}}}{A_{\infty}} , (1+2^{-k}) \frac{ a_{i} \ell_{n_{k+1}}}{A_{\infty}}\right]
\end{equation} 
simultaneously for all $k \geq 1$ and all branches $ \mathsf{b}_{i} \in \mathcal{B}_{k}$.
 \end{corollary} 
 
 Note that this implies what we have said previously: if the sequence $(n_k)_{k\geq 1}$ grows sufficiently fast, then the event $\{ \mathcal{B}_{k} \ne \varnothing, \forall k \geq 1\}$ has a probability at least $1-\eps$.

\proof This is a direct application of Lemma \ref{lemma:tech}. Imagine that $n_{1}, \ldots , n_{k}$ have been fixed and that $ \mathcal{B}_{k}$ has been constructed and is non empty. Fix $ \mathsf{b}_{i} \in \mathcal{B}_{k}$. Using Lemma \ref{lemma:tech} (i) with $S =  \mathsf{b}_{i}$, $n = n_{k+1}$  and $\eta = 2^{-k}$,  we get 
  \begin{eqnarray*} \mathbb{P}\left( \left| | \mathcal{T}( \mathsf{b}_{i}) \cap \mathcal{B}_{k+1}| - \frac{ a_{i} \ell_{n_{k+1}}}{A_{\infty}}\right| \geq 2^{-k} \frac{ a_{i} \ell_{n_{k+1}}}{A_{\infty}} \right)   &\leq& 4^{k} (n_{k+1})^{-c + \circ(1)}/a_{i}\\  & \underset{ \mathsf{b}_{i} \mathrm{\ is \ good}}{\leq}&  4^{k} (n_{k+1})^{-c + \circ(1)}n_{k}^{\alpha + \varepsilon}.  \end{eqnarray*}
Given $n_{k}$, we can thus choose $f(n_{k})$ large enough so that if $n_{k+1} \geq f(n_{k})$ the right-hand side of the last display is at most $ 2^{-k} \varepsilon /(n_{k}+1)$. For such an integer $n_{k+1}$, the probability that one of the branches $ \mathsf{b}_{i}$ of $ \mathcal{B}_{k}$ does not satisfy (\ref{eq:taille}) is at most 
$$ (n_{k}+1) \cdot 2^{-k} \varepsilon /(n_{k}+1) = 2^{-k} \varepsilon.$$ Constructing in this way a sequence $(n_{k})_{k\geq 1}$, we see that  the probability that (\ref{eq:taille}) fails for one $k$ is at most $ \varepsilon \cdot (2^{-1}+ 2^{-2}+...) = \varepsilon$.\endproof 

\begin{lemma} \label{lem:nombre}
There exists a function $g : \mathbb{N} \to \mathbb{N}$ with $g(n)> 2n$  for all $n \geq 1$, such that if the sequence $(n_k)_{k \geq 1}$ satisfies $n_{k+1} \geq g(n_{k})$ for all $k \geq 1$, then with probability at least $1- \varepsilon$,  for all $k \geq 1$ we have 
$$ \sup_{x \in \mathcal{T}} \#\left\{ \mathsf b_i \in \mathcal B_k : \mathsf{b}_{i} \cap B(x,n_{k}^{-\alpha}) \ne \varnothing \right\} \leq n_{k}^{ \varepsilon}$$

\end{lemma}
\proof Imagine that $ \mathcal{B}_{k}$ is constructed and pick $\mathsf b_{i} \in \mathcal{B}_{k}$. Conditionally on the number $N$ of branches of $ \mathcal{B}_{k+1}$ grafted onto $ \mathsf{b}_{i}$, the grafting points of these branches are i.i.d.~and uniform on $ \mathsf{b}_{i}$. We decompose the good branch $ \mathsf{b}_{i}$ into $  \lceil {a_{i}}/{n_{k+1}^{-\alpha}} \rceil$ intervals of length at most $n_{k+1}^{-\alpha}$. If none of these intervals contains more than $n_{k+1}^{\varepsilon/2}$ branches then it is not possible to have more than $3 n_{k+1}^{\varepsilon/2}$ branches within distance less than $n_{k+1}^{-\alpha}$. Noticing that $N \leq n_{k+1}+1$, we get that the probability to have more than $3 n_{k+1}^ {\varepsilon/2}$ branches within distance less than $n_{k+1}^{-\alpha}$ is at most
\begin{eqnarray*}
\left\lceil \frac{a_{i}}{n_{k+1}^{-\alpha}} \right \rceil \cdot {N \choose n_{k+1}^{\varepsilon/2}} \left( \frac{n_{k+1}^{-\alpha}}{a_{i}} \right)^{n_{k+1}^{\varepsilon/2}}  &\leq&   \left( \frac{n_{k}^{-\alpha + \circ(1)}}{n_{k+1}^{-\alpha}} +1\right) \cdot (n_{k+1}+1)^{n_{k+1}^ {\varepsilon/2}} n_{k+1}^{- \alpha  \cdot n_{k+1}^{\varepsilon/2}} n_{k}^{ (\alpha + \varepsilon) \cdot n_{k+1}^{ \varepsilon/2}} \\ &\leq & \big(n_{k+1}^{-\alpha + 1 + \circ(1)}  n_{k}^{ \alpha + \varepsilon + \circ(1)}\big)^{n_{k+1}^{\varepsilon/2}}. 
 \end{eqnarray*}
Clearly by making $n_{k+1} \geq g(n_{k})$ grows rapidly enough we can ensure that the series of the last probabilities is as small as we wish. Hence with probability at least $1- \varepsilon$, for every $k \geq 2$ and any $x \in \mathcal{T}$, the number of branches of $ \mathcal{B}_{k}$ grafted on a given $ \mathsf{b}_{i} \in \mathcal{B}_{k-1}$ within distance $ n_{k}^{-\alpha}$ of $x$ is at most $3n_{k}^{\varepsilon/2}$. Using this proposition in cascades (and remarking that $n_{i}^{-\alpha} > n_{k}^{-\alpha}$ for $i < k$), we get that on this event 
 $$ \sup_{x \in \mathcal{T}} \#\left\{ \mathsf b_i \in \mathcal B_k : \mathsf{b}_{i} \cap B(x,n_{k}^{-\alpha}) \ne \varnothing \right\} \leq  3n_{1}^{\varepsilon/2}  \cdots 3n_{k-1}^{   \varepsilon/2} 3 n_{k}^{ \varepsilon/2},$$  
and the last product is at most $n_{k}^{\varepsilon}$ provided that $n_{k}$ grows rapidly enough.
\endproof 

We will now use this lemma and Lemma \ref{lemma:tech} to control the maximal length of groups of branches of $\mathcal B_{k+1}$ that are grafted on a ball of radius  $r$, when the center of the ball runs over $\mathcal B_k$. In that aim, we also need to assume that the sequence $(n_k)$ grows sufficiently fast so that
\begin{equation}
\label{eq:encoreplusgros}
n_{k} = n_{k+1}^{\circ(1)} \quad \mbox{ as } k \to \infty. 
\end{equation}

\begin{corollary} 
\label{cor:size}
Assume that the sequence $(n_k)$ satisfies $n_{k+1} \geq g(n_k)$ for all $k$ -- where $g$ is the function of the previous lemma -- as well as (\ref{eq:encoreplusgros}).
For each $k \in \mathbb N$, each  $r>0$ and each $x \in \mathcal B_k$, consider the total length of branches  of $\mathcal B_{k+1}$ that are grafted on $B(x,r) \cap \mathcal B_k \subset \mathcal T_{n_{k+1}-1}$. Let $\mathcal L_{k+1}(r)$ be the supremum of these lengths when $x$ runs over $\mathcal B_k$. Then with probability at least $1-\eps$, for all  $0<\gamma<1-\eps/\alpha$ and for all $k$ large enough (the threshold depending on $\gamma)$,
$$
\mathcal L_{k+1}(r) \leq r^{\frac{1}{\alpha}-\eps} \ell_{n_{k+1}} \quad \text{for all } r \in \Big[n_{k+1}^{-\alpha}, n_{k+1}^{-1+\frac{\eps}{2}}\Big]
$$
and
\begin{equation*}
\label{ineggrandr}
\mathcal L_{k+1}(r) \leq r^{\gamma}\ell_{n_{k+1}} \quad \text{for all } r \in \Big[n_{k+1}^{-1+\frac{\eps}{2}}, n_k^{-\alpha}\Big].
\end{equation*} 
\end{corollary}

\proof
Let $\mathcal A$ denote the event of probability at least $1-\eps$ on which the conclusion of Lemma \ref{lem:nombre} holds. In the following, we will work mostly on $\mathcal A$ and $\gamma \in (0,1-\eps/\alpha)$ is fixed.

To start with, we set up for each $r \in [n_{k+1}^{-\alpha},n_k^{-\alpha}]$ a specific covering of $\mathcal B_k$. 
Split each $\mathsf b_i \in \mathcal B_{k}$ into $\lceil a_i/r\rceil$ intervals, with $\lfloor a_i/r\rfloor$ intervals of length $r$ and a last one (if $a_i/r$ is not an integer) of length at most $r$ which is chosen to be the one that reaches the leaf of $\mathsf b_i$. This gives a set of $$\sum_{i:\mathsf b_i \in \mathcal B_k} \left\lceil \frac{a_i}{r}\right\rceil \ \leq \ \frac{|\mathcal B_k|}{r} + \# G_{n_k} \ \leq \ \frac{A_{\infty}}{r}  + n_k+1$$ intervals of $\mathcal B_k$ of lengths at most $r$. Besides, consider the balls of radius $r$ centered at the points of $\mathcal B_{k-1} \cap \mathcal B_k$ (i.e. at the ``roots" of the $\mathsf b_i, \mathsf b_i\in \mathcal B_k$). For such a ball $B$, the set $B \cap \mathcal B_k$ intersects at most $n_k^{\eps}$ branches $\mathsf b_i, \mathsf b_i\in \mathcal B_k$, conditionally on $\mathcal A$ (by Lemma \ref{lem:nombre}). In particular, its length $|B \cap \mathcal B_k|$ is at most $n_k^{\eps}r$. The covering we are interested in is composed by the intersections of these balls with $\mathcal B_k$ and the intervals mentioned above. It is therefore composed by sets that all have a length at most  $n_k^{\eps}r$. Moreover, each ball of radius $r$ centered at a point of $\mathcal B_k$ is included in the union of  two neighboring elements of the covering, one of which being necessarily an interval.
\medskip

\noindent $\bullet$ Using this covering, we note that 
\begin{eqnarray*}
&&\mathbb P\left(\exists r \in \Big[n_{k+1}^{-\alpha}, n_{k+1}^{-1+\frac{\eps}{2}}\Big] : \mathcal L_{k+1}(r) \geq r^{\frac{1}{\alpha}-\eps} \ell_{n_{k+1}}, \mathcal A \right) \\
&\leq& \mathbb P\left( \mathcal L_{k+1}\big(n_{k+1}^{-1+\frac{\eps}{2}}\big) \geq (n_{k+1}^{-\alpha})^{\frac{1}{\alpha}-\eps} \ell_{n_{k+1}}, \mathcal A \right) \\ 
&\leq&  \left(A_{\infty}n_{k+1}^{1-\frac{\eps}{2}} +n_{k}+1\right) \cdot 2\mathbb P\left( \mathcal X \geq 2^{-1}n_{k+1}^{-1+\alpha \eps} \ell_{n_{k+1}} \right),
\end{eqnarray*}
where $\mathcal X$ represents the total length of branches  of $\mathcal B_{k+1}$ that are grafted on a subset $S \subset \mathcal T_{n_{k+1}-1}$ of length $n_k^{\eps}n_{k+1}^{-1+\eps/2}$. By Lemma \ref{lemma:tech} (ii), for all integers $m\geq 1$ and then all $k$ large enough, we have  
\begin{eqnarray*}
\mathbb P\left(\exists r \in \Big[n_{k+1}^{-\alpha}, n_{k+1}^{-1+\frac{\eps}{2}}\Big] : \mathcal L_{k+1}(r) \geq r^{\frac{1}{\alpha}-\eps} \ell_{n_{k+1}}, \mathcal A \right) &\leq& C'_m \left(A_{\infty}n_{k+1}^{1-\frac{\eps}{2}} +n_{k}+1\right) \frac{n_{k}^{\eps m} n_{k+1}^{(-1+\eps/2)m} \ell_{n_{k+1}}^m}{n_{k+1}^{(-1+\alpha\eps)m} \ell_{n_{k+1}}^m} \\
&\leq & n_{k+1}^{1-\frac{\eps}{2}+(\frac{1}{2}-\alpha) \eps m + \circ(1)}.
\end{eqnarray*}
Fix $m$ large enough so that  the exponent $1-\eps/2+(1/2-\alpha) \eps m \leq -1$. Since  $n_{k+1} \geq 2^k$ for all $k$, we can therefore use Borel-Cantelli's lemma to conclude that  on $\mathcal A$, almost surely for  all $k$ large enough,
$$
\mathcal L_{k+1}(r) \leq r^{\frac{1}{\alpha}-\eps} \ell_{n_{k+1}} \quad \text{for all } r \in \Big[n_{k+1}^{-\alpha}, n_{k+1}^{-1+\frac{\eps}{2}}\Big].
$$

\noindent $\bullet$ For  $r \in [n_{k+1}^{-1+\eps/2},n_k^{-\alpha}]$ the argument is similar but we have to split the interval $[n_{k+1}^{-1+\eps/2},n_k^{-\alpha}]$ into subintervals to conclude. Let $\eta \in (1,(1-\eps \alpha^{-1})/\gamma)$ and first note that
\begin{eqnarray*}
\mathbb P\left(\exists r \in \Big[n_{k+1}^{-1+\frac{\eps}{2}},n_k^{-\alpha}\Big] : \mathcal L_{k+1}(r) \geq r^{\gamma} \ell_{n_{k+1}}, \mathcal A\right) &\leq & \sum_{n=0}^{N_k} \mathbb P\left( \exists r \in \Big[n_k^{-\alpha\eta^{n+1}},n_k^{-\alpha \eta^n}\Big] : \mathcal L_{k+1}(r) \geq r^{\gamma} \ell_{n_{k+1}}, \mathcal A\right) \\
&\leq & \sum_{n=0}^{N_k} \mathbb P\left(\mathcal L_{k+1}(n_k^{-\alpha \eta ^n}) \geq n_k^{-\alpha \gamma \eta^{n+1}}  \ell_{n_{k+1}}, \mathcal  A\right),
\end{eqnarray*}
where $N_k$ is the largest integer $n$ such that $n_k^{-\alpha \eta^n} \geq n_{k+1}^{-1+\eps/2}$. Applying Lemma \ref{lemma:tech} (ii) to subsets $S$ of $\mathcal T_{n_{k+1}-1}$ of lengths $n_k^{\eps} n_k^{-\alpha\eta^n}$, we see that for all integers $m\geq 1$ and then all $k$ large enough and all $n \leq N_k$,  \begin{eqnarray*}
\mathbb P\left(\mathcal L_{k+1}(n_k^{-\alpha\eta^n}) \geq n_k^{-\alpha\gamma\eta^{n+1}} \ell_{n_{k+1}}, A\right) &\leq& C_m \left(A_{\infty}n_k^{\alpha\eta^n} +n_k+1\right) \frac{\big(n_k^{\eps}n_k^{-\alpha\eta^n}\big)^m \ell_{n_{k+1}}^m}{\big(n_k^{-\alpha\gamma \eta^{n+1}}\big)^m \ell_{n_{k+1}}^m} 
\\
&\leq &  C'_m n_k^{\left(\alpha+(\eps+\alpha(\gamma \eta-1))m\right)\eta^n},
\end{eqnarray*}
where we have used for the last inequality that $\eta^n \geq 1$ and $\alpha>1$. 
The parameters have been chosen so that $\eps+\alpha(\gamma \eta-1)<0$. So we can fix $m$ sufficiently large so that $\alpha+(\eps+\alpha(\gamma \eta-1))m \leq -1$ and then conclude that for all $k$ large enough
\begin{eqnarray*}
\mathbb P\left(\exists r \in \Big[n_{k+1}^{-1+\frac{\eps}{2}},n_k^{-\alpha}\Big] : \mathcal L_{k+1}(r) \geq r^{\gamma} \ell_{n_{k+1}}, A\right) &\leq & C'_m \sum_{n=0}^{N_k} \frac{1}{n_k^{\eta^n}} \\
&\underset{n_k \geq 2^{k-1}}{\leq} &  \frac{C'_m}{2^{(k-1)}} \sum_{n=0}^{\infty} \frac{1}{2^{(k-1)(\eta^n-1)}} \\
& \leq &   \frac{C'_m}{2^{(k-1)}} \sum_{n=0}^{\infty} \frac{1}{2^{\eta^n-1}}
\end{eqnarray*}
and the series, clearly, is convergent. Again, we conclude with Borel-Cantelli's lemma that a.s. on $\mathcal A$, for all $k$ large enough, 
$$
\mathcal L_{k+1}(r) \leq r^{\gamma}\ell_{n_{k+1}} \quad \text{for all } r \in \Big[n_{k+1}^{-1+\frac{\eps}{2}}, n_k^{-\alpha}\Big].$$
\endproof

\subsubsection{Proof of Proposition \ref{prop:lowbound}} Fix $\gamma \in (1-\eps,1-\eps/\alpha)$ and fix a sequence $(n_{k})_{k\geq1}$ such that the conditions of Corollary \ref{cor:main} and Corollary \ref{cor:size} are satisfied (in particular (\ref{eq:encoreplusgros}) holds). There exists therefore an event $ \mathcal{E}$ of probability at least $1- 2\varepsilon$ on which the conclusions of Lemma \ref{lem:leaves}, Corollary \ref{cor:main} and Corollary \ref{cor:size} hold, for the $\gamma$ we have chosen. From now on, we work on this event $\mathcal E$
and it is implicit in what follows that all assertions hold conditionally on $\mathcal E$.
By Corollary \ref{cor:main}, each branch of $ \mathcal{B}_{k}$ will have some branches of $ \mathcal{B}_{k+1}$ grafted on it and so $ \mathcal{B}_{k} \ne \varnothing$ for all $k \geq 1$ and the measures $ \pi_{k}$ are well-defined for all $k \geq 1$.  We denote by $\pi$ a subsequential limit of $(\pi_{k})$. We aim at proving \eqref{eq:lowerboundgoal}.

By Corollary \ref{cor:main} again, for all $k \geq 1$
 \begin{eqnarray} | \mathcal{B}_{k+1}| \in \left[(1-2^{-k}) \frac{ | \mathcal{B}_{k}| \ell_{n_{k+1}}}{A_{\infty}} , (1+2^{-k}) \frac{ | \mathcal{B}_{k}| \ell_{n_{k+1}}}{A_{\infty}}\right].   \label{eq:rapide2}\end{eqnarray}
Consequently, 
\begin{eqnarray} \label{eq:rapide} |\mathcal B_{k+1}| \underset{ \eqref{eq:goodbranches}}{=} n_{k+1}^{1-\alpha+\circ(1)}| \mathcal{B}_{k}| \underset{ \eqref{eq:encoreplusgros}}{=} n_{k+1}^{1-\alpha+\circ(1)}.  \end{eqnarray}
Next, using Corollary \ref{cor:main} as well as  \eqref{eq:rapide2} in cascades, we see that for any $ \mathsf{b}_{i} \in \mathcal{B}_{k}$ and any $k' \geq k$
 \begin{eqnarray*} | \mathcal{T}( \mathsf{b}_{i}) \cap \mathcal{B}_{k'}| &\in& a_{i} \cdot \left[ \prod_{j=k+1}^{k'} (1-2^{-(j-1)})\frac{\ell_{n_{j}}}{A_{\infty}} ; \prod_{j=k+1}^{k'} (1+2^{-(j-1)})\frac{\ell_{n_{j}}}{A_{\infty}} \right]. \\
 |\mathcal{B}_{k'}| &\in& | \mathcal{B}_{k}|\cdot  \left[ \prod_{j=k+1}^{k'} (1-2^{-(j-1)})\frac{\ell_{n_{j}}}{A_{\infty}} ; \prod_{j=k+1}^{k'} (1+2^{-(j-1)})\frac{\ell_{n_{j}}}{A_{\infty}} \right]. \end{eqnarray*}
Let $c_{1} = \prod_{j=1}^\infty (1-2^{-j})/(1+2^{-j}) \in (0,\infty)$ and $c_{2} = \prod_{j=1}^\infty (1+2^{-j})/(1-2^{-j}) \in (0,\infty)$, then we have 
   $$ \pi_{k'}( \mathcal{T}( \mathsf{b}_{i})) = \frac{| \mathcal{T}( \mathsf{b}_{i}) \cap \mathcal{B}_{k'}|}{ | \mathcal{B}_{k'}|}  \in \frac{a_{i}}{| \mathcal{B}_{k}|} \cdot [ c_{1},c_{2}].$$
 Using arguments similar to those developed in the proof of Lemma \ref{lem:leaves} we get that for any branch $ \mathsf{b}_{i} \in \mathcal{B}_{k}$ 
 \begin{eqnarray} \label{eq:pimasse} \pi( \mathcal{T}( \mathsf{b}_{i}))  \in \left[ \frac{c_{1}}{c_{2}} \frac{a_{i}}{| \mathcal{B}_{k}|},  \frac{c_{2}}{c_{1}} \frac{a_{i}}{| \mathcal{B}_{k}|}\right]. 
 \end{eqnarray} 
Now, recall that the support of the measure $\pi$ is included in $\cap_{k\geq 1} \mathcal T(\mathcal B_k)$ (by Lemma \ref{lem:leaves}) and fix $x \in \cap_{k\geq 1} \mathcal T(\mathcal B_k)$. Let $r \in [n_{k+1}^{-\alpha},n_k^{-\alpha}]$ for some $k \in \mathbb N$ and note that 
 \begin{eqnarray*}
 \pi \left(B(x,r) \right)&=&\sum_{i:\mathsf b_i \in \mathcal B_{k+1}}  \pi \left(B(x,r)\cap \mathcal T(\mathsf b_i) \right) \\
 &\underset{(\ref{eq:pimasse})}\leq & \frac{c_2}{c_1 | \mathcal{B}_{k+1}|}\sum_{i:\mathsf b_i \in \mathcal B_{k+1}}  a_i \mathbbm 1_{\{B(x,r)\cap \mathcal T(\mathsf b_i) \neq \varnothing\}}.
 \end{eqnarray*}
Note also that $\sum_{i:\mathsf b_i \in \mathcal B_{k+1}}  a_i \mathbbm 1_{\{B(x,r)\cap \mathcal T(\mathsf b_i) \neq \varnothing\}} \leq \mathcal L_{k+1}(r)$, with the notation of Corollary \ref{cor:size}. (The bounds below will therefore be true simultaneously for all $x$.) Hence, according to this corollary, 
 $$
 \pi(B(x,r)) \leq   \frac{c_2 \ell_{n_{k+1}} r^{\frac{1}{\alpha}-\eps}}{c_1 | \mathcal{B}_{k+1}|} \underset{(\ref{eq:rapide2})} \leq  \frac{c_2 A_{\infty}}{c_1 (1-2^{-k})} \cdot \frac{r^{\frac{1}{\alpha}-\eps}}{| \mathcal{B}_{k}| }    \leq r^{1/\alpha-3\eps/2} \quad \text{ for all }r \in \Big[n_{k+1}^{-\alpha},n_{k+1}^{-1+\frac{\eps}{2}}\Big]
 $$
 provided that $k$ is large enough, since $| \mathcal{B}_{k}|=n_k^{1-\alpha+\circ(1)}=n_{k+1}^{\circ(1)}$, by (\ref{eq:rapide}) and (\ref{eq:encoreplusgros}). On the other hand, again by Corollary \ref{cor:size}, 
$$
 \pi \left(B(x,r) \right) \leq  \frac{c_2 A_{\infty}}{c_1 (1-2^{-k})} \cdot \frac{r^{\gamma}}{| \mathcal{B}_{k}| }   \underset{(\ref{eq:rapide})}=  \frac{r^{\gamma}}{n_k^{1-\alpha+\circ(1)}} \quad \text{for all }r \in \Big[n_{k+1}^{-1+\frac{\eps}{2}}, n_{k}^{-\alpha}\Big],
 $$
 where the $\circ(1)$ is independent of $r$. Recall that $\gamma>1-\eps$ and then note that $r \leq n_k^{-\alpha}$ implies $r^{\gamma-1/\alpha+\eps} \leq n_{k}^{1-\alpha \gamma -\alpha \eps}$, hence $r^{\gamma}n_k^{-1+\alpha+\circ(1)}\leq r^{1/\alpha-\eps}$ for all $k$ large enough (independently of $r \leq n_k^{-\alpha}$).

 In conclusion, on the event $\mathcal E$, for all $k$ large enough and then all $r \in [n_{k+1}^{-\alpha},n_k^{-\alpha}]$ -- hence for all $r$ sufficiently small, $$\pi(B(x,r))\leq r^{1/\alpha-3\eps/2} \quad  \text{ for all } x \in \bigcap_{k\geq 1} \mathcal T(\mathcal B_k),$$  which  implies \eqref{eq:lowerboundgoal} since the support of $ \pi$ is included in $\cap_{k\geq 1} \mathcal T(\mathcal B_k)$. 

\section{Appendix} \label{appendix}

We gather here some elementary technical results useful in the core of the paper.
Let $(a_i,i\geq 1)$ be a sequence of strictly positive real numbers, and $A_i=a_1+\ldots+a_i$, $i\geq 1$.

\begin{lemma} \label{lem:manip} 
Assume that $0<a_i \leq c$ for all $i\geq 1$ and some $c<\infty$. Then, 
\begin{enumerate} 
\item[$\mathrm{(i)}$] the series $\sum_{i} \frac{a_i}{A^2_i}$ and $\sum_{i} \big(\frac{a_i}{A_i}\big)^2$ are convergent
\item[$\mathrm{(ii)}$] if $a_i \leq i^{-\alpha+\circ(1)}$ for some $\alpha>0$, then $\sum_{i\geq n} \frac{a_i^2}{A_i} \leq n^{-\alpha+\circ(1)}$
\smallskip
\item[$\mathrm{(iii)}$] if $A_i \geq i^{1-\alpha+\circ(1)}$ for some $\alpha \in (0,1)$, then $\sum_{i\geq n} \frac{a_i}{A_i^2} \leq n^{\alpha-1+\circ(1)}$.
\end{enumerate}
\end{lemma}

\proof
Since the sequence $(A_i^{-1})$ is bounded from above, Assertions (i) and (ii) are immediate when the series $\sum_i a_i$ is convergent. (Assertion (iii) requires anyway that the series $\sum_i a_i$ is divergent.) 
So \text{we assume from now on that the series $\sum_i a_i$ diverges}, and define for all $k\geq 1$
$$
n_k:=\inf\{i\geq 1 : A_i \geq k\},
$$
which is finite.
Note that $A_{n_k}\geq k$ and $A_{n_{k+1}-1} <k+1$, in particular $A_{n_{k+1}-1}-A_{n_k} <1$ and therefore $\sum_{i=n_k}^{n_{k+1}-1} a_i <c+1$.

\medskip

\noindent \textit{Assertion} (i). The convergence of the series $\sum_i \frac{a_i}{A_i^2}$ is simply due to the following observation : 
\begin{eqnarray*}
\sum_{i=n_1}^{\infty}\frac{a_i}{A_i^2} \ = \ \sum_{k=1}^{\infty} \sum_{i=n_k}^{n_{k+1}-1} \frac{a_i}{A_i^2} \ \leq \ \sum_{k=1}^{\infty} \frac{1}{k^2}  \sum_{i=n_k}^{n_{k+1}-1} a_i \ < \ \sum_{k=1}^{\infty} \frac{c+1}{k^2}.
\end{eqnarray*}
The convergence of the series $\sum_i \big(\frac{a_i}{A_i}\big)^2$ follows, since $\big(\frac{a_i}{A_i}\big)^2 \leq  \frac{ca_i}{A_i^2}.$
\medskip

\noindent \textit{Assertion} (ii). We assume that $a_i \leq i^{-\alpha+\circ(1)}$ for some $\alpha \in (0,1]$. Let $\eps \in (0,\alpha/2)$. For $i$ large enough, we have $A_i \leq  i^{1-\alpha+\eps}$ and therefore, for $k$ large enough, $n_k \geq k^{1/(1-\alpha+\eps)}$. Consequently, for all $i \geq \max(n,n_k)$, with $n$ and $k$ large enough, 
$$
a_i \ \leq \ i^{-\alpha+\eps} \ = \ i^{-\alpha + 2\eps} \times i^{-\eps} \ \leq  n^{-\alpha + 2\eps} \times k^{-\eps/(1-\alpha-\eps)}.
$$
And then, for $n$ large enough,
\begin{eqnarray*}
\sum_{i\geq n} \frac{a_i^2}{A_i}&=&\sum_{k \geq 1} \sum_{i=n_k}^{n_{k+1}-1} \mathbbm 1_{\{i\geq n\}}\frac{a_i^2}{A_i} \\
&\leq& n^{-\alpha+2\eps} \sum_{k \geq 1} \frac{k^{-\eps/(1-\alpha+\eps)}}{k} \left( \sum_{i=n_k}^{n_{k+1}-1} a_i \right) \\
&\leq& n^{-\alpha+2\eps} \sum_{k \geq 1} \frac{c+1}{k^{1+\eps/(1-\alpha+\eps)}}. 
\end{eqnarray*}
This holds for all $\eps>0$ small enough and the conclusion follows. 

\medskip

\noindent \textit{Assertion} (iii). Fix $\eps \in (0,(1-\alpha)/2)$. For $i$ large enough, $A_i \geq i^{1-\alpha-\eps}$.  Hence for $i \geq \max(n,n_k)$, with $n$ large enough, 
$$
A_i^2 \ \geq \ A_n^{1-\eps} A_{n_k}^{1+\eps} \ \geq \ n^{1-\alpha-2\eps} k^{1+\eps}.
$$
Consequently, for $n$ large enough
$$
\sum_{i \geq n} \frac{a_i}{A_i^2} \ = \ \sum_{k=1}^{\infty} \sum_{i=n_k}^{n_{k+1}-1} \mathbbm 1_{\{i\geq n\}}  \frac{a_i}{A_i^2} \ \leq \ n^{\alpha-1+2\eps} \sum_{k=1}^{\infty} \frac{c+1}{k^{1+\eps}}.$$
\endproof

\bibliographystyle{siam}
\bibliography{bibli.bib}

\end{document}